\documentclass[12pt,a4paper]{article}
\usepackage{graphicx}
\usepackage{amsfonts}
\usepackage{float}
\usepackage[latin1]{inputenc}
\usepackage{blkarray}
\usepackage{color}
\usepackage{amsmath}
\usepackage{amstext}
\usepackage{amssymb}
\usepackage{array}
\usepackage{multicol}
\usepackage{pdfpages}
\usepackage{multirow}
\usepackage{threeparttable}

\usepackage{subcaption}
\def\N{\Bbb N}
\def\R{\Bbb R}

\def\D{\cal D}

\def\<{\langle}
\def\>{\rangle}

\def\Chi{\raise .3ex \hbox{\large $\chi$}}

\def\n{{\bf }}

\def\[{\Bigl [}
  \def\]{\Bigr ]}
\def\({\Bigl (}
\def\){\Bigr )}
\def\[{\Bigl [}
  \def\]{\Bigr ]}
\def\({\Bigl (}
\def\){\Bigr )}

\def\div{{\mbox{\rm div}}}
\def\dsp{\displaystyle}
\def\K{{\bf K}}
\def\x{{\bf x}}

\def\n{{\bf n}}
\def\s{{\bf s}}

\def\g{{\bf g}}

\def\V{{\bf V}}

\def\q{{\bf q}}

\def\G{{\Gamma}}

\def\D{{\cal D}}
\def\cells{{\cal M}}
\def\faces{{\cal F}}
\def\nodes{{\cal V}}
\def\edges{{\cal E}}

\allowdisplaybreaks

\begin{document}


  \title{Parallel numerical modeling of hybrid-dimensional compositional non-isothermal Darcy flows in fractured porous media}
  
  \author{
F. Xing\thanks{Laboratoire de Math\'ematiques J.A. Dieudonn\'e, UMR 7351 CNRS, University Nice Sophia Antipolis, team COFFEE, INRIA Sophia Antipolis M\'editerran\'ee, Parc Valrose 06108 Nice Cedex 02, France, and BRGM Orl\'eans France, feng.xing@unice.fr},         
R. Masson\thanks{Laboratoire de Math\'ematiques J.A. Dieudonn\'e, UMR 7351 CNRS, University Nice Sophia Antipolis, and team COFFEE, INRIA Sophia Antipolis M\'editerran\'ee, Parc Valrose 06108 Nice Cedex 02, France, roland.masson@unice.fr}, 
S. Lopez\thanks{BRGM, scientific and Technical Center, 3 avenue Claude Guillemin, BP 36009, 45060 Orl\'eans Cedex 2 France, s.lopez@brgm.fr}
}
\date{}

\maketitle

%
%
%
%
  

  \begin{abstract}
This paper introduces a new discrete fracture model accounting for non-isothermal 
compositional multiphase Darcy flows and 
complex networks of fractures with intersecting, immersed and non immersed fractures. 
The so called hybrid-dimensional  model using a 2D model in the fractures 
coupled with a 3D model in the matrix is first derived rigorously starting from the equi-dimensional matrix fracture 
model. Then, it is discretized using a fully implicit time integration combined with the 
Vertex Approximate Gradient (VAG) finite volume scheme which is adapted to polyhedral meshes 
and anisotropic heterogeneous media. The fully coupled systems are assembled and solved in parallel 
using the Single Program Multiple Data (SPMD) paradigm with one layer of ghost cells. This strategy allows for a 
local assembly of the discrete systems. 
An efficient preconditioner is implemented to solve the linear systems at each time step and 
each Newton type iteration of the simulation. The numerical efficiency of our approach is assessed on 
different meshes, fracture networks, and physical settings in terms of parallel scalability, 
nonlinear convergence and linear convergence.  

  \end{abstract}




\section{Introduction}
\label{sec_intro}

Flow and transport in fractured porous media are of paramount importance for many applications such as 
petroleum exploration and production, geological storage of carbon dioxide, hydrogeology, or geothermal energy. 
Two classes of models, dual continuum and discrete fracture models, 
are typically employed and possibly coupled to simulate flow and transport in fractured porous media. 
Dual continuum models assume that the fracture network is well connected  and can 
be homogeneized as a continuum coupled  
to the matrix continuum using transfer functions. On the other hand, 
discrete fracture models (DFM), on which this paper focuses, represent explicitly the fractures as codimension one 
surfaces immersed in the surrounding matrix domain. 
The use of lower dimensional rather than equi-dimensional entities to represent the fractures has been introduced in 
\cite{MAE02}, \cite{FNFM03}, \cite{BMTA03}, \cite{KDA04}, \cite{MJE05} 
to facilitate the grid generation and to reduce the number of degrees of freedom of the discretized model.  
The reduction of dimension in the fracture network is obtained from the equi-dimensional model by 
integration and averaging along the width of each fracture. 
The resulting so called hybrid-dimensional {models} couple the 3D model in the matrix with a 2D model in the fracture network 
taking into account the jump of the normal fluxes as well as additional transmission conditions  at the matrix fracture interfaces. 
These transmission conditions 
depend on the mathematical nature of the equi-dimensional model and on additional physical assumptions.   
They are typically derived for a single phase Darcy flow for which they 
specify either the continuity of the 
pressure in the case of fractures acting as drains (see \cite{MAE02}, \cite{GSDFN}) or  Robin type 
conditions in order to take into account the discontinuity of the pressure for fractures acting as barriers 
(see \cite{FNFM03}, \cite{MJE05}, \cite{ABH09}, \cite{BHMS2016}). 
Different transmission conditions are derived in \cite{tracer2016} in the case 
of a linear hyperbolic equation, and in \cite{RJBH06}, \cite{HF08}, \cite{Jaffre11}, \cite{BGGM14}, \cite{BHMS16}
in the case of two-phase immiscible Darcy flows. \\

The discretization of hybrid-dimensional Darcy flow models  
has been the object of many  works. In \cite{KDA04} a cell-centred Finite Volume scheme using a 
Two Point Flux Approximation (TPFA) is proposed assuming the orthogonality of the mesh and 
isotropic permeability fields. Cell-centred Finite Volume schemes 
can be extended to general meshes and anisotropic permeability fields 
using  MultiPoint Flux Approximations (MPFA) 
(see \cite{HADEH09}, \cite{TFGCH12}, \cite{SBN12}, \cite{AELHP152D}, \cite{AELHP153D}). 
MPFA schemes can lack robustness on distorted meshes and large anisotropies due to 
the non symmetry of the discretization. They are also very expensive compared with 
nodal discretizations on tetrahedral meshes.
In \cite{MAE02}, a Mixed Finite Element (MFE) method is proposed for single phase 
Darcy flows. It is extended to two-phase flows in \cite{HF08} in an IMPES framework using 
a Mixed Hybrid Finite Element (MHFE) discretization 
for the pressure equation and a Discontinuous Galerkin discretization of the saturation equation.
The Hybrid Finite Volume and Mimetic finite difference schemes, belonging  
to the family of Hybrid Mimetic Mixed Methods (HMM) \cite{DEGH13},  
have been extended to hybrid-dimensional models in \cite{FFJR16}, 
\cite{AFSVV16} as well as in \cite{GSDFN}, \cite{BHMS2016} 
in the more general Gradient Discretization framework \cite{DEGGH16}. 
These approaches are adapted to general meshes and anisotropy 
but require as many degrees of freedom as faces.  
Control Volume Finite Element Methods (CVFE) 
\cite{BMTA03}, \cite{RJBH06}, \cite{MF07}, \cite{HADEH09} have the advantage to 
use only nodal unknowns leading to much fewer degrees of freedom than 
MPFA and HMM schemes on tetrahedral meshes. On the other hand, at the matrix fracture interfaces, 
the control volumes 
have the drawback to be shared between the matrix and the fractures. It results that a 
strong refinement of the mesh is needed at these interfaces in the case of large contrasts 
between the matrix and fracture permeabilities. 
This article focus on the Vertex Approximate Gradient (VAG) scheme which 
has been introduced for the discretization of multiphase Darcy flows 
in \cite{EHGM-CG-12} and extended to hybrid-dimensional models in 
\cite{BGGM14}, \cite{GSDFN}, \cite{BHMS2016}, \cite{tracer2016}, \cite{BHMS16}.  
The VAG scheme uses nodal and fracture face unknowns 
in addition to the cell unknowns which can be eliminated without any fill-in. 
Thanks to its essentially nodal feature, it leads to a sparse discretization on tetrahedral or mainly tetrahedral meshes. 
It has the advantage, compared 
with the CVFE methods  of \cite{BMTA03}, \cite{RJBH06}, \cite{MF07} or \cite{MMB2007}, 
to avoid the mixing of the control volumes at the  matrix fracture interfaces, 
which is a key feature for its coupling with a transport model. As shown 
in \cite{BGGM14} for two-phase flow problems, this allows 
for a coarser mesh size at the matrix fracture interface for a given accuracy.
Let us also mention that non-matching discretizations of the fracture and matrix meshes
are studied for single phase Darcy flows in \cite{DS12}, \cite{FS13}, \cite{BPS14} 
and \cite{SFHW15}. \\

The first objective of this paper is 
to extend the derivation of the hybrid-dimensional model to the case of 
non-isothermal compositional multiphase Darcy flows. 
{To focus on compositional non-isothermal features, capillary pressures are not considered 
in this paper. They could be included following the usual phase based upwinding approach as in \cite{EHGM-CG-12} or recent ideas {developed} in \cite{BGJMP16} for two-phase flows. 
Let us refer to \cite{BHMS17} for a comparison of both approaches in the case of an immiscible two-phase flow using a reference solution provided by the equi-dimensional model in the fractures.} All the underlying assumptions of our reduced model will be carefully stated. 
In particular, the fractures are considered as pervious and are assumed not to act as barriers. 
It results, as in \cite{MAE02}, that the pressure can be considered as continuous at the matrix fracture interfaces. 
The hybrid-dimensional model accounts for complex network of fractures including 
intersecting, immersed and non immersed fractures. 
The formulation of the compositional model is based on a Coats' type  
formulation \cite{coats-89-imp}, \cite{CHB02} extending the approach presented in \cite{EHGM-CG-12} to non-isothermal flows. 
It accounts for an arbitrary nonzero number of components in each phase allowing to 
model immiscible, partially miscible or fully miscible flows. 

The second objective of this paper is to extend the VAG discretization to our 
model and to develop an efficient 
parallel algorithm implementing the discrete model. 
Following \cite{EHGM-CG-12}, \cite{BGGM14}, the discretization is based on a finite volume formulation 
of the component molar and energy conservation equations. The definition of the control volumes is adapted 
to the heterogeneities of the porous medium and avoids in particular the mixing of matrix and fracture 
rocktypes for the degrees of freedom located at the matrix fracture interfaces. 
The fluxes combine the VAG Darcy and Fourier fluxes with  
a phase based upwind approximation of the mobilities. A fully implicit Euler 
time integration coupling the conservation equations with the local closure laws including thermodynamical equilibrium 
is used in order to avoid severe limitations on the time step due 
to the high velocities and small control volumes in the fractures. 

The discrete model is implemented in parallel based 
on the SPMD (Single Program, Multiple Data) paradigm. 
It relies on a distribution of the mesh on the processes with one layer of ghost cells in order 
to allow for a local assembly of the discrete systems. The key ingredient for the efficiency of the parallel 
algorithm is the solution, at each time step and at each Newton type iteration, of the large sparse linear 
system coupling the physical unknowns on the spatial degrees of freedom of the VAG scheme. Our strategy is first based 
on the elimination, without any fill-in, of both the local closure laws and the cell unknowns. 
Then, the reduced linear system is solved using a parallel iterative solver preconditioned by a 
CPR-AMG preconditioner introduced in \cite{LVW:2001} and \cite{SMW:2003}. This state of the art 
preconditioner  combines multiplicatively an Algebraic MultiGrid (AMG) preconditioner for a proper pressure block of the linear system with 
a local incomplete factorization preconditioner for the full system. 
The numerical efficiency of the algorithm, in terms of parallel scalability, nonlinear convergence 
and linear convergence, 
is investigated on several test cases. We consider 
different families of meshes and different complexity of fracture networks ranging 
from a few fractures to say about 1000 fractures 
with highly contrasted matrix fracture permeabilities. 
The test cases incorporate different physical models including one isothermal immiscible two-phase flow, one  
isothermal Black Oil two-phase flow model, as well as three non-isothermal water component liquid gas flow models.  \\

This paper is organized as follows. In section \ref{sec_model}, 
the  hybrid-dimensional  non-isothermal compositional multiphase Darcy flow model is derived from the equi-dimensional model. 
In Section \ref{sec_algo}, the VAG discretization is briefly recalled and then 
extended to our model. The parallel algorithm is detailed in section \ref{sec_pl}. 
Section \ref{sec_num} is devoted to the test cases including the numerical investigation of the parallel scalability of the algorithm.

\section{Hybrid-dimensional compositional non-isothermal Darcy flow model}
\label{sec_model}

This section deals with the modeling and the formulation of non-isothermal compositional multiphase Darcy 
flows in fractured porous media. 
The fractures are represented as surfaces of co-dimension one immersed in the surrounding three dimensional  
matrix domain. The 2D Darcy flow in the fracture network is coupled with the 3D Darcy flow 
in the matrix domain, hence the terminology of hybrid-dimensional model. The reduction of dimension 
in the fracture is obtained by extension to non-isothermal compositional flows 
of the methodology introduced in \cite{MAE02}, \cite{FNFM03}, \cite{MJE05} 
for single phase Darcy flows. Complex networks of fractures are considered 
including immersed, non immersed and intersecting planar fractures.  
The formulation of the 
compositional model is based on a Coats' type 
formulation \cite{coats-89-imp} extending to non-isothermal flows 
the approach presented in \cite{EHGM-CG-12}. 
It accounts for an arbitrary nonzero number of components in each phase allowing to 
model immiscible, partially miscible or fully miscible flows. 
{To focus on the compositional and non-isothermal aspects, 
we consider a Darcy flow model without capillary pressures. The capillary pressures including 
different rocktypes at the matrix fracture interface can be taken 
into account in the framework of the VAG scheme following the usual phase based upwinding 
of the mobilities as in \cite{EHGM-CG-12}. An alternative approach is proposed 
in \cite{BGJMP16} for two-phase flows in order to capture the jump of the saturations at the matrix fracture interface $\Gamma$ 
due to discontinuous capillary pressure curves.
These two choices are compared in \cite{BHMS17} to a reference solution provided by an equi-dimensional model in the fractures.
It is shown that the second choice provides a better solution as long as the matrix acts as a barrier since it captures the saturation jump. On the other hand, the first choice provides a more accurate solution when the non wetting phase goes out of the fractures since
the mean capillary pressure in the fractures is better approximated.}

\subsection{Extended Coats' formulation of non-isothermal compositional models}
\label{subsec_coats}

Let us denote by ${\cal P}$ the set of phases and by ${\cal C}$ the set of components. 
Each phase $\alpha\in {\cal P}$ 
is described by its non empty subset of components 
${\cal C}^\alpha\subset {\cal C}$ in the sense that it contains 
the components $i\in {\cal C}^\alpha$. 
It is assumed that, for any $i\in {\cal C}$, the set of phases containing the component $i$
$$
{\cal P}_i = \{\alpha\in {\cal P}\,|\, i\in {\cal C}^\alpha\}.
$$ 
is non empty. 
The thermodynamical properties of each phase $\alpha \in {\cal P}$  depend on the pressure $P$, the temperature $T$, and the molar fractions
$$
C^\alpha = \left(C_i^\alpha\right)_{i\in {\cal C}^\alpha}. 
$$
For each phase $\alpha \in {\cal P}$, we denote 
by $\zeta^\alpha(P,T,C^\alpha)$ the molar density, by  
 $\rho^\alpha(P,T,C^\alpha)$ the mass density, by 
 $\mu^\alpha(P,T,C^\alpha)$ the dynamic viscosity, by 
 $f_i^\alpha(P,T,C^\alpha)$, $i\in {\cal C}^\alpha$ the fugacity coefficients, by 
 $e^\alpha(P,T,C^\alpha)$ the molar internal energy, and by 
 $h^\alpha(P,T,C^\alpha)$ the molar enthalpy. The relative permeabilities are denoted 
for each phase $\alpha\in {\cal P}$ by $k_r^\alpha(S)$ where $S=(S^\alpha)_{\alpha\in {\cal P}}$ is 
the vector of the phase volume fractions (saturations). 
{The model takes into account phase change reactions 
which are assumed to be at equilibrium. It results that phases can appear or disappear. }
Therefore, we denote by $Q \subset {\cal P}$, $Q\neq \emptyset$ 
the unknown representing 
the set of present phases. 
For a given set of present phases $Q$, it may occur that 
a component $i\in{\cal C}$ does not belong to the subset $\bigcup_{\alpha \in Q} {\cal C}^\alpha$
of ${\cal C}$. 
Hence, we define the subset of absent components as a function of $Q$ by 
\begin{equation*}
\overline {\cal C}_Q = \{i\in{\cal C}\,|\, Q\cap{\cal P}_i = \emptyset\}. 
\end{equation*}
Following \cite{coats-89-imp}, \cite{CHB02}, \cite{EHGM-CG-12}, 
the extended non-isothermal Coats' formulation relies on the the so-called natural variables and 
uses the set of unknowns
$$
X = \(P,T,S^\alpha, C^\alpha, \alpha\in Q, n_i, i\in \overline {\cal C}_Q, Q\).  
$$
The saturations are implicitely set to $S^\alpha = 0$ for all absent phases $\alpha\in {\cal P}\setminus Q$. 
Let us denote by $n_i(X)$ the number of moles of the component $i\in {\cal C}$ 
per unit pore volume defined as the independent unknown $n_i$ for $i\in \overline{\cal C}_Q$ and as 
\begin{equation*}
n_i(X) = \sum_{\alpha\in Q \cap {\cal P}_i} \zeta^{^\alpha}(P,T, C^{\alpha})~S^{\alpha}~C_i^{\alpha} 
\end{equation*}
for $i\in {\cal C}\setminus\overline{\cal C}_Q$.
The fluid energy per unit pore volume is denoted by 
$$
E(X) = \sum_{\alpha \in Q} \zeta^\alpha(P,T,C^\alpha) S^\alpha e^\alpha(P,T,C^\alpha),
$$
and the rock energy per unit rock volume is denoted by $E_r(P,T)$. 
For each phase $\alpha \in Q$, we denote by $m_i^\alpha(X)$ the mobility 
of the component $i\in {\cal C}^\alpha$ in phase $\alpha\in Q$ with 
$$
m_i^\alpha(X) = C^\alpha_i \zeta^\alpha(P,T,C^\alpha)  \frac{k_r^\alpha(S)}{\mu^\alpha(P,T,C^\alpha) }, 
$$
and by 
$$
m^\alpha_e(X)= h^\alpha(P,T,C^\alpha) \zeta^\alpha(P,T,C^\alpha)  
\frac{k_r^\alpha(S)}{\mu^\alpha(P,T,C^\alpha) }
$$
the flowing enthalpy in phase $\alpha\in Q$. 
The generalized Darcy velocity of the phase $\alpha\in Q$ is
$$
\frac{k_r^\alpha(S)}{\mu^\alpha(P,T,C^\alpha) } {\bf V}^\alpha
\ \mathrm{with} \
{\bf V}^\alpha = - {\bf K} \( \nabla P - \rho^\alpha(P,T,C^\alpha) {\bf g}\), 
$$
where ${\bf g}$ is the gravitational acceleration.  
The total molar flux of the component $i\in {\cal C}\setminus \overline{\cal C}_Q$ is denoted by 
$$
{\bf q}_i = \sum_{\alpha\in Q \cap {\cal P}_i} m^\alpha_i(X) {\bf V}^\alpha, 
$$
and the energy flux is defined as 
$$
{\bf q}_e = \sum_{\alpha \in Q} m_e^\alpha(X) {\bf V}^\alpha - \lambda \nabla T, 
$$ 
where $\lambda$ is the thermal conductivity of the fluid and rock mixture. 

The system of equations accounts for the molar conservation for each component $i \in {\cal C}$ and  
the energy conservation 
\begin{equation}
\label{cons_coats}
\begin{split}
\phi \partial_t n_i & + \mathrm{div}({\bf q}_i) = 0, \, i \in {\cal C},\\
\phi \partial_t E + (1-\phi) \partial_t E_r & +\div ( {\bf q}_e ) = 0,
\end{split}
\end{equation}
coupled to the following 
local closure laws including the thermodynamical equilibrium for each component $i$ 
present in at least two phases among the set of present phases  $Q$
\begin{equation}
\label{closure_coats}
\begin{split}
\dsp{ \sum_{\alpha\in Q} S^{\alpha} } = & 1,  \\
\dsp{\sum_{i\in {\cal C}^\alpha} C_i^{\alpha}} = & 1, \,\, \alpha \in Q,\\
f_i^\alpha(P,T,C^\alpha)C_i^\alpha = & f_i^\beta(P,T,C^\beta)C_i^\beta, \,\,  
\alpha\neq\beta,  (\alpha,\beta)\in (Q\cap{\cal P}_i)^2. 
\,\,i\in{\cal C}.
\end{split}
\end{equation}
The system is closed with an additional equation for the discrete unknown  $Q$ 
which is typically  obtained by a flash calculation or by simpler criteria depending on the specific 
thermodynamical system. It provides the fixed point equation denoted by 
$$
Q = Q_{flash}(X). 
$$
\subsection{Discrete fracture network}
\label{subsec_dfn}

Let $\Omega$ denote a bounded domain of $\R^3$  assumed to be polyhedral. 
Following \cite{MAE02}, \cite{FNFM03}, \cite{MJE05}, \cite{GSDFN}, \cite{BHMS2016} the fractures are represented as interfaces of codimension 1.  
Let $J$ be a finite set and let 
$
\overline \Gamma = \bigcup_{j\in J} \overline \Gamma_j
$  
and its interior $\Gamma = \overline \Gamma\setminus \partial\overline\Gamma$ 
denote the network of fractures $\Gamma_j\subset \Omega$, $j\in J$, such that each $\Gamma_j$ is 
a planar polygonal simply connected open domain included in an oriented plane of $\R^3$. 
\begin{figure}[htbp!]
\begin{center}
\includegraphics[width=0.37\textwidth]{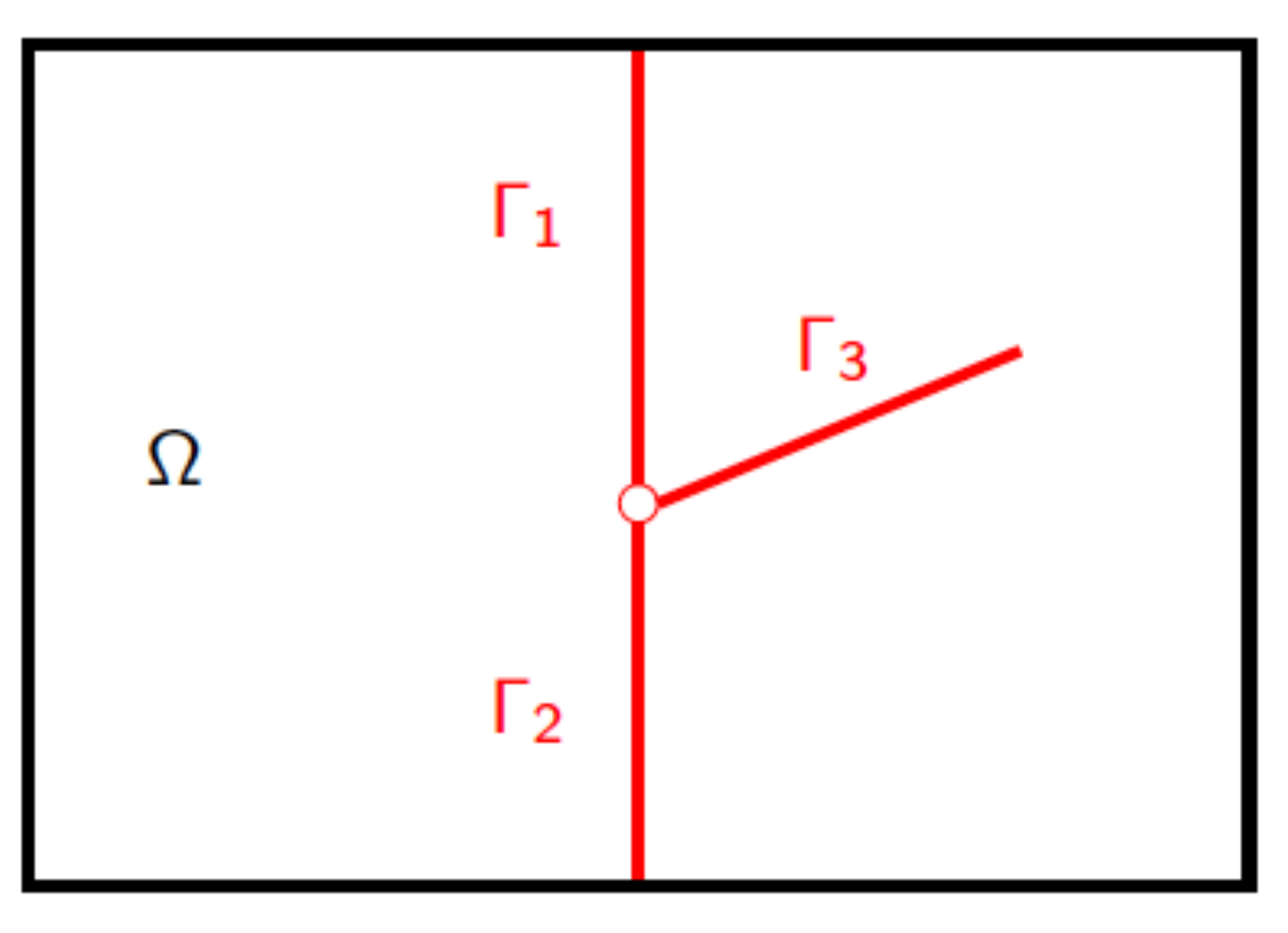}   
\caption{Example of a 2D domain with 3 intersecting fractures $\Gamma_1, \Gamma_2, \Gamma_3$.}
\label{fig_network}
\end{center}
\end{figure}
The fracture width is denoted by $d_f$ and is such that $0 < \underline d_f \leq d_f(\x) \leq \overline d_f$ 
for all $\x\in \Gamma$. 
We can define, for each fracture $j\in J$, its two 
sides $+$ and $-$. For scalar functions on $\Omega$, possibly 
discontinuous at the interface $\Gamma$ (typically in $H^1(\Omega\setminus \overline\Gamma)$, 
we denote by $\gamma^\pm$ the trace operators on the side $\pm$ of $\Gamma$. 
Continuous scalar functions $u$ at the interface $\Gamma$ (typically in $H^1(\Omega)$) 
are such that $\gamma^+ u = \gamma^- u$ 
and we denote by $\gamma$ the trace operator on $\Gamma$ for such functions. 
At almost every point of the fracture network, 
we denote by ${\bf n}^\pm$ the unit normal vector oriented outward to the side $\pm$ of $\Gamma$ 
such that ${\bf n}^+ + {\bf n}^- = 0$. For vector fields on $\Omega$,  possibly discontinuous at the interface $\Gamma$ (typically in $H_\div(\Omega\setminus\overline\Gamma)$, 
we denote by $\gamma_n^{\pm}$ the normal trace operator on the side $\pm$ of $\Gamma$ 
oriented outward to the side $\pm$ of $\Gamma$. 

The gradient operator in the matrix domain $\Omega\setminus \overline\Gamma$ is denoted by $\nabla$ and the tangential gradient operator on the 
fracture network is denoted by  $\nabla_\tau$ such that 
$$
\nabla_\tau u  = \nabla u - (\nabla u\cdot \n^+)\n^+. 
$$
We also denote by $\div_{\tau}$ the tangential divergence operator on the fracture network, and by 
$d\tau({\bf x})$ the Lebesgue measure on $\Gamma$.

We denote by $\Sigma$ the dimension $1$ open set defined by the intersection of the fractures 
excluding the boundary of the domain $\Omega$, i.e. the interior of 
$\bigcup_{\{(j,j')\in J\times J\,|\, j\neq j'\}} \partial \Gamma_j \cap \partial \Gamma_{j'} 
\setminus \partial \Omega$.  

For the matrix domain, Dirichlet (subscript $D$) and Neumann (subscript $N$) boundary conditions are imposed on 
the two dimensional open sets $\partial \Omega_D$ and $\partial \Omega_N$ respectively where 
$\partial \Omega_D\cap \partial \Omega_N = \emptyset$, 
$\partial \Omega = \overline {\partial \Omega_D}\cup \overline {\partial \Omega_N}$.
Similarly for the fracture network, the Dirichlet and Neumann boundary conditions are imposed on the one dimensional open sets $\partial \Gamma_D$ and $\partial \Gamma_N$ respectively where 
$\partial \Gamma_D\cap \partial \Gamma_N = \emptyset$, 
$\overline{\partial \Gamma \cap \partial\Omega} = \overline {\partial \Gamma_D}\cup \overline {\partial \Gamma_N}$.

\subsection{Hybrid-dimensional  model for a two-phase flow example}
\label{subsec_hybrid_twophase}

For the sake of clarity, in this subsection, we first extend the hybrid-dimensional model proposed in
\cite{MAE02},  \cite{FNFM03}, \cite{MJE05} 
in the case of a single phase Darcy flow to the case of a two-phase Darcy flow. 

Let us denote by $\<\>_\Gamma$ the averaging operator in the width of the fracture in the normal direction, 
and let us set $S^\alpha_f = \<S^\alpha\>_\Gamma$ and $P_f = \<P\>_\Gamma$. 
In order to explain our construction of the 
hybrid-dimensional model from the equi-dimensional model, let us consider the following  
immiscible, incompressible isothermal two-phase flow model 
on the full domain $\Omega$ with 3D fractures 
$$
\phi \partial_t S^\alpha + \div( \q^\alpha ) = 0, \, \alpha =1,2,
$$
with constant dynamic viscosities $\mu^\alpha$, $\alpha=1,2$ and the following Darcy two-phase velocities 
$$
\q^\alpha = - { k_r^\alpha(S^\alpha) \over \mu^\alpha } {\bf K} (\nabla P -\rho^\alpha \g ), \, \alpha =1,2,  
$$
where $S^1 + S^2 = 1$. 
The permeability tensor ${\bf K}$ is assumed to be constant in the width of the fractures 
and to have the normal vector $\n^+$ as principal direction. We denote by ${\bf K}_f$ 
the corresponding tangential permeability tensor and by $k_{f,n}$ the corresponding 
normal permeability, both defined as function of $\x\in \Gamma$.  
The porosity is also assumed to be constant in the width of the fracture and 
denoted by $\phi_f$ as a function defined on $\Gamma$.

The reduction of dimension in the fractures is based on the assumption that 
$d_f \ll \mbox{diam}(\Omega)$. It is obtained by integration of 
the conservation equations along the width of the fractures in the normal direction 
using the approximation $k_r^\alpha(S^\alpha_f)$ of 
$k_r^\alpha(S^\alpha)$  in the definition of the tangential flux in the fractures  
\begin{equation}
\label{red_eqfrac}
\begin{split}
&d_f \phi_f \partial_t S^\alpha_f + \div_\tau (\q^\alpha_f) 
- \gamma_n^+ \q^\alpha_m -\gamma_n^- \q^\alpha_m  = 0, \, \alpha=1,2,\\
&\q^\alpha_f = {k_r^\alpha(S_f^\alpha) \over \mu^\alpha }\V^\alpha_f, \, \alpha=1,2,\\
&\V^\alpha_f = - d_f {\bf K}_f ( \nabla_\tau P_f -\rho^\alpha \g_\tau ), \, \alpha=1,2,\\
& S^1_f + S^2_f = 1,
\end{split}
\end{equation}
with $\g_\tau = \g - (\g\cdot\n^+)\n^+$. 
This conservation equation on $\Gamma$ is coupled to the conservation equation 
in the matrix domain $\Omega\setminus \overline \Gamma$  
\begin{equation}
\label{red_eqmat}
\begin{split}
&\phi_m \partial_t S_m^\alpha + \div( \q^\alpha_m ) = 0, \, \alpha=1,2,\\
& \q_m =  {k_r^\alpha(S_m^\alpha) \over \mu^\alpha } \V^\alpha_m, \, \alpha=1,2,\\
& \V^\alpha_m = - {\bf K}_m (\nabla P_m -\rho^\alpha \g ), \, \alpha=1,2,\\
& S^1_m + S^2_m = 1,
\end{split}
\end{equation}
where we use the subscript $m$ to denote the properties and unknowns of the reduced model 
defined in the matrix domain.  

This hybrid-dimensional model \eqref{red_eqfrac}-\eqref{red_eqmat} 
is closed with transmissions 
conditions at the matrix fracture interface $\Gamma$. They are based on the following two-point 
flux approximation of the normal fluxes  at both sides $\pm$ 
of the fractures. As opposed to the model proposed in \cite{Jaffre11}, 
we take into account the gravity term which cannot be neglected. 
\begin{equation}
\label{red_tpfa}
V^{\alpha,\pm}_{f,n} = {2 k_{f,n} \over d_f} (\gamma^\pm P_m - P_f) 
+ k_{f,n} \rho^\alpha  \g\cdot \n^\pm,   \, \alpha=1,2. 
\end{equation}
{In  \cite{Jaffre11}, the definition of the normal fluxes $\gamma_n^\pm \q^\alpha_m$ 
is obtained with the mobility ${k_r^\alpha(S^\alpha_f) \over \mu^\alpha }$ using the mean saturation 
in the width of the fracture. This choice  cannot account for
the propagation of the saturation front from the matrix to the fracture.
To solve this problem, we propose to use a monotone two point flux between the interface on the matrix side and the centre of the fracture. Our choice is based on the phase based upwind flux leading to upwind the mobility with respect to  
the sign of   $V^{\alpha,\pm}_{f,n}$. }
For any $a\in \R$ let us set $a^+ = \max(a,0)$ and $a^- = \min(a,0)$. 
The normal fluxes are obtained using the following upwind approximations 
of the mobilities with respect to the sign of the phase normal velocities
\begin{equation}
\label{red_upwind}
\gamma_n^\pm \q^\alpha_m = {k_r^\alpha(S^\alpha_f) \over \mu^\alpha } (V^{\alpha,\pm}_{f,n})^- 
+   { k_r^\alpha(\gamma^\pm S^\alpha_m) \over \mu^\alpha }  (V^{\alpha,\pm}_{f,n})^+, 
\, \alpha=1,2. 
\end{equation}
This phase based upwinding of the mobilities is known to lead to a two point monotone flux for the saturation equation. It also provides a flux consistency error of the order of 
the ratio between the width of the fracture and the size of the matrix domain, which is assumed 
to be small.  

Note also that the use of the mean saturation $S^\alpha_f$ 
in the mobilities \eqref{red_upwind}  for output fluxes from the fracture to the matrix basically assumes 
that the saturation in the fracture is well approximated by a constant along the width. This holds 
true for fractures with a high conductivity $d_f {\bf K}_f$ compared with 
the conductivity of the matrix $diam(\Omega) {\bf K}_m$. This condition will be assumed in the following. \\

Moreover, when ${k_{f,n}\over d_f} \gg {K_m \over diam(\Omega)}$, the transmission condition \eqref{red_tpfa} 
can be further approximated by the pressure continuity condition at the matrix 
fracture interface $\Gamma$
\begin{equation}
\label{red_continuousP}
\gamma^+ P_m = \gamma^- P_m = \gamma P_m = P_f,   
\end{equation}
recovering the condition introduced in \cite{MAE02} for single phase Darcy flows. In this case, the definition 
of the normal fluxes \eqref{red_upwind} is modified as follows using the normal trace  $\gamma_n^\pm \V^\alpha_m$ 
of $\V^\alpha_m$ rather than $V_{f,n}^{\alpha,\pm}$:  
\begin{equation}
\label{red_continuous_upwind}
\gamma_n^\pm \q^\alpha_m = {k_r^\alpha(S^\alpha_f) \over \mu^\alpha } (\gamma_n^\pm \V^\alpha_m)^- 
+   { k_r^\alpha(\gamma^\pm S^\alpha_m) \over \mu^\alpha }  (\gamma_n^\pm \V_m^\alpha)^+, 
\, \alpha=1,2. 
\end{equation}
In the following, we will assume that this approximation holds which means that we consider 
the case of fractures acting as drains and exclude
the case of \textcolor{red}{fractures} acting as barriers. \\

Finally, closure conditions are set at the immersed boundary of the fracture network 
$\partial \Gamma \setminus \partial \Omega$ (fracture tips) as well as 
at the intersection $\Sigma$ between fractures. 
Let $\gamma_{n_{\partial \Gamma}}$ (resp. $\gamma_{n_{\partial \Gamma_j}}, j\in J$) 
denote the normal trace operator at the fracture network boundary (resp. fracture $\Gamma_j$ boundary) 
oriented outward to $\Gamma$ (resp. $\Gamma_j$). 
At fracture tips, it is classical to assume homogeneous Neumann boundary conditions in the sense that  
$$
\gamma_{n_{\partial \Gamma}} \q^\alpha_f = 0, \mbox{ on } \partial \Gamma \setminus \partial \Omega, \, \alpha=1,2,
$$ 
meaning that the flow at the tip of a fracture 
can be neglected compared with the flow along the sides of the fracture. 
At the fracture intersection $\Sigma$, we introduce the additional unknowns 
$P_\Sigma$, $S_\Sigma^\alpha$,  $\alpha=1,2$ and we impose 
the normal flux conservation equations  
$$
\sum_{j\in J} (\gamma_{n_{\partial \Gamma_j}} \q^\alpha_{f} )|_\Sigma =0, \,\, \alpha=1,2, 
$$
meaning that the volume at the intersection between fractures is neglected. The saturations 
$S_\Sigma^\alpha$, $\alpha=1,2$ are such that $S_\Sigma^1 + S_\Sigma^2 = 1$ and play the role of 
the input saturations at the fracture intersection. 
In addition, we also impose the continuity of the pressure $P_f = P_\Sigma$ at $\Sigma$. This amounts to assume 
a high ratio between the permeability at the intersection and the fracture width 
compared with the ratio between the tangential permeability of each fracture and its lengh. 

\subsection{Hybrid-dimensional non-isothermal compositional model}
\label{subsec_hybrid_comp}

The hybrid-dimensional non-isothermal compositional model is obtained following 
the above strategy for the dimension reduction.
The set of unknowns is defined by $X_m$ in the matrix domain $\Omega\setminus\overline\Gamma$, 
by $X_f$ in the fracture network $\Gamma$, and by $X_\Sigma$ at the fracture intersection $\Sigma$. 
The set of equations 
couples the molar and energy conservation equations in the matrix 
\begin{equation}
\label{hybrid_cons_mat}
\begin{split}
\phi_m \partial_t n_i(X_m)  + \div({\bf q}_{i,m})  =&  0, \, i\in {\cal C},\\
\phi_m \partial_t E(X_m) + (1-\phi_m) \partial_t E_r(P_m,T_m)  
+\div ( {\bf q}_{e,m} )  =& 0,
\end{split}
\end{equation}
in the fracture network 
\begin{equation}
\label{hybrid_cons_frac}
\begin{split}
d_f \phi_f \partial_t n_i(X_f)  + & \div_\tau({\bf q}_{i,f}) - \gamma_n^+ \q_{i,m} - \gamma_n^- \q_{i,m}  = 0, \, i\in {\cal C},\\
d_f \phi_f \partial_t E(X_f) + & d_f (1-\phi_f) \partial_t E_r(P_f,T_f)   \\
+ & \div_\tau ( {\bf q}_{e,f} ) - \gamma_n^+ \q_{e,m} - \gamma_n^- \q_{e,m} = 0,
\end{split}
\end{equation}
and at the fracture intersection
\begin{equation}
\label{hybrid_cons_sigma}
\begin{split}
\sum_{j\in J} (\gamma_{n_{\partial \Gamma_j}} \q_{i,f})|_{\Sigma} =&0, \, i\in {\cal C},\\
\sum_{j\in J} (\gamma_{n_{\partial \Gamma_j}} \q_{e,f})|_{\Sigma} =&0,  
\end{split}
\end{equation}
as well as the Darcy and Fourier laws providing the fluxes in the matrix 
\begin{equation}
\label{Darcy_Fourier_mat}
\begin{split}
\q_{i,m} =&\dsp \sum_{\alpha\in Q_m \cap {\cal P}_i} m^\alpha_i(X_m) {\bf V}_m^\alpha, \\
\q_{e,m} = & \sum_{\alpha \in Q_m} m_e^\alpha(X_m) {\bf V}_m^\alpha 
- \lambda_m \nabla T_m, 
\end{split}
\end{equation}
and in the fracture network 
\begin{equation}
\label{Darcy_Fourier_frac}
\begin{split}
\q_{i,f} =&\dsp \sum_{\alpha\in Q_f \cap {\cal P}_i} m^\alpha_i(X_f) {\bf V}_f^\alpha, \\
\q_{e,f} = & \sum_{\alpha \in Q_f} m_e^\alpha(X_f) {\bf V}_f^\alpha 
- d_f \lambda_f \nabla_\tau T_f, 
\end{split}
\end{equation}
where
\begin{equation*}
\begin{split}
{\bf V}^\alpha_m & = - {\bf K}_m \( \nabla P_m - \rho^\alpha(P_m,T_m,C_m^\alpha) {\bf g}\),\\
{\bf V}^\alpha_f & = - d_f {\bf K}_f \( \nabla_\tau P_f - \rho^\alpha(P_f,T_f,C_f^\alpha) {\bf g}_\tau\), 
\end{split}
\end{equation*}
and finally the local closure laws including the thermodynamical equilibrium 
\begin{equation}
\label{Hybrid_closure}
\begin{split}
\dsp{ \sum_{\alpha\in Q_\nu} S^{\alpha}_\nu } = & 1,  \\
\dsp{\sum_{i\in {\cal C}^\alpha} C_{i,\nu}^{\alpha}} = & 1, \,\, \alpha \in Q_\nu,\\
f_i^\alpha(P_\nu,T_\nu,C_\nu^\alpha)C_{i,\nu}^\alpha = & f_i^\beta(P_\nu,T_\nu,C_\nu^\beta)C_{i,\nu}^\beta, \,\,  
\alpha\neq\beta,  (\alpha,\beta)\in (Q_\nu\cap{\cal P}_i)^2. 
\,\,i\in{\cal C},\\
Q_\nu = & Q_{flash}(X_\nu), 
\end{split}
\end{equation}
for $\nu=m,f,\Sigma$. The system \eqref{hybrid_cons_mat}-\eqref{hybrid_cons_frac}-\eqref{hybrid_cons_sigma}-\eqref{Darcy_Fourier_mat}-\eqref{Darcy_Fourier_frac}-\eqref{Hybrid_closure} 
is closed with the transmission conditions at the matrix fracture interface $\Gamma$.  
These conditions state, as above, the continuity of the pressure complemented for non-isothermal 
models with the continuity  of the temperature. 
It is combined with a phase based upwind approximation of the mobilities 
in the matrix fracture normal fluxes. This corresponds to the usual 
finite volume two point upwind scheme for the mobilities (see e.g. \cite{aziz-settari-79}) 
applied for our reduced model in the normal direction 
between the centre of the fracture and each side of the fracture.  
\begin{equation}
\label{Hybrid_transmission}
\begin{split}
& \gamma^+ P_m= \gamma^- P_m = \gamma P_m = P_f,\\
& \gamma^+ T_m= \gamma^- T_m = \gamma T_m = T_f,\\
& \gamma_n^\pm \q_{i,m} = m_i^\alpha(X_f) (\gamma_n^\pm \V^\alpha_m)^- 
+ m_i^\alpha(\gamma^\pm X_m) (\gamma_n^\pm \V_m^\alpha)^+,\\
& \gamma_n^\pm \q_{e,m} = m_e^\alpha(X_f) (\gamma_n^\pm \V^\alpha_m)^- 
+ m_e^\alpha(\gamma^\pm X_m) (\gamma_n^\pm \V_m^\alpha)^+ + \gamma_n^\pm (-\lambda_m \nabla T_m). 
\end{split}
\end{equation}
Note also that the pressure $P_f$ (resp. the temperature $T_f$)  
is assumed continuous and equal to $P_\Sigma$ (resp. $T_\Sigma$) 
at the fracture intersection $\Sigma$, and that 
homogeneous Neumann boundary conditions are applied 
for each component molar $\q_{i,f}$ and energy $\q_{e,f}$ 
fluxes at the fracture tips $\partial \Gamma\setminus \Omega$. 

Regarding the boundary conditions, to fix ideas, we restrict ourselves
to either Dirichlet or homogeneous Neumann boundary conditions. 
At the Dirichlet matrix boundary $\partial \Omega_D$ (resp. Dirichlet fracture boundary $\partial \Gamma_D$) 
the pressure $P_{m,D}$ (resp. $P_{f,D}$), temperature $T_{m,D}$ (resp. $T_{f,D}$),
are specified, as well as the set of input phases $Q_{m,D}$ (resp. $Q_{f,D}$), their 
volume fractions $S^\alpha_{m,D}$, $\alpha\in Q_{m,D}$ (resp.  $S^\alpha_{f,D}$, $\alpha\in Q_{f,D}$) 
and their molar fractions  $C^\alpha_{m,D}$, $\alpha\in Q_{m,D}$ (resp.  $C^\alpha_{f,D}$, $\alpha\in Q_{f,D}$) 
assumed to satisfy the local closure laws.  
Then, we set for $\nu=m,f$ 
\begin{equation}
\label{Hybrid_transmission_Cont}
\begin{split}
 P_{\nu} =& P_{\nu,D},\\
 T_{\nu} =& T_{\nu,D},\\
 S^\alpha_\nu =&  S^\alpha_{\nu,D} \mbox{ for } \alpha\in Q_{\nu,D} \mbox{ if } \V_{\nu}^\alpha\cdot \n_\nu < 0,\\
 C^\alpha_\nu =&  C^\alpha_{\nu,D} \mbox{ for } \alpha\in Q_{\nu,D} \mbox{ if } \V_{\nu}^\alpha\cdot \n_\nu < 0, 
\end{split}
\end{equation}
where $\n_\nu$ is the output unit normal vector at the boundary $\partial \Omega_D$ for $\nu=m$, 
and at the boundary $\partial \Gamma_D$ for $\nu = f$. 

Homogeneous Neumann boundary conditions are applied at the boundaries 
$\partial \Omega_N$ and $\partial \Gamma_N$ in the sense 
that $\q_{i,\nu}\cdot \n_\nu = 0$  for $i\in {\cal C}\cup \{e\}$, $\nu=m,f$,  
where $\n_\nu$ is the output unit normal vector at the boundary $\partial \Omega_N$ for $\nu=m$, and at the boundary $\partial \Gamma_N$ for $\nu = f$.

\section{Discretization and algorithm}
\label{sec_algo}

\subsection{VAG discretization}

The VAG discretization of hybrid-dimensional two-phase Darcy flows introduced in \cite{BGGM14} considers generalised polyhedral meshes of $\Omega$ in the spirit of  \cite{Eymard.Herbin.ea:2010}. {In short, the mesh is assumed conforming, the cells are star-shapped polyhedrons, and faces are not necessarily planar in the sense that they
can be defined as the union of triangles joining the edges of the face to a so-called face centre. 
In more details,} let $\cells$ be the set of cells that are disjoint open polyhedral subsets of $\Omega$ such that
$\bigcup_{K\in\cells} \overline{K} = \overline\Omega$, for all $K\in\cells$, ${\x}_K$ denotes the so-called ``centre'' of the cell $K$ under the assumption that $K$ is star-shaped with respect to ${\x}_K$. 
The set of faces of the mesh is denoted by $\faces$ and $\faces_K$ is the set of faces of the cell $K\in \cells$. 
The set of edges of the mesh is denoted by $\edges$ and $\edges_\sigma$ is the set of edges of the face $\sigma\in \faces$. 
The set of vertices of the mesh is denoted by $\nodes$ and $\nodes_\sigma$ is the set of vertices of the face $\sigma$. 
For each $K\in \cells$ we define $\nodes_K = \bigcup_{\sigma\in \faces_K} \nodes_\sigma$.   

The faces are not necessarily planar. It is just assumed that for each face $\sigma\in\faces$, there exists a so-called ``centre'' of the face ${{\bf x}}_\sigma \in {\sigma}\setminus \bigcup_{e\in \edges_\sigma} e$ such that
$
{\x}_\sigma = \sum_{\s\in \nodes_\sigma} \beta_{\sigma,\s}~\x_\s, \mbox{ with }
\sum_{\s\in \nodes_\sigma} \beta_{\sigma,\s}=1,
$
and $\beta_{\sigma,\s}\geq 0$ for all $\s\in \nodes_\sigma$; moreover
the face $\sigma$ is assumed to be defined  by the union of the triangles
$T_{\sigma,e}$ defined by the face centre ${\x}_\sigma$
and each edge $e\in\edges_\sigma$. 
The mesh is also supposed to be conforming w.r.t. the fracture network $\G$ in the sense that for all  $j\in J$ there exist the subsets $\faces_{\G_j}$ of $\faces$ such that 
$$
\overline \G_j = \bigcup_{\sigma\in\faces_{\G_j}} \overline{\sigma}.  
$$
We will denote by $\faces_\G$ the set of fracture faces 
$$
\faces_\Gamma = \bigcup_{j\in J} \faces_{\G_j}, 
$$ 
and by 
$$
\nodes_\Gamma = \bigcup_{\sigma\in \faces_\Gamma} \nodes_{\sigma}, 
$$ 
the set of fracture nodes. 
This geometrical discretization of $\Omega$ and $\G$ is denoted in the following by $\D$. 

In addition, the following notations will be used 
$$
\cells_\s = \{K\in \cells\,|\, \s\in \nodes_K\}, 
\
\cells_{\sigma} = \{K\in \cells\,|\, \sigma\in \faces_K\}, 
$$
and 
$$
\faces_{\G,\s} = \{\sigma \in \faces_\Gamma \,|\, \s\in \nodes_\sigma\}. 
$$

\begin{figure}[H]
\begin{center}
\includegraphics[width=0.4\textwidth]{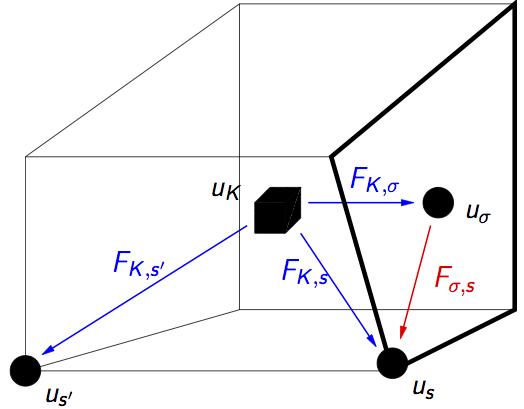}  
\caption{For a cell $K$ and a fracture face $\sigma$ (in bold), 
examples of VAG degrees of freedom $u_K$, $u_\s$, $u_\sigma$, $u_{\s'}$ and 
VAG fluxes $F_{K,\sigma}$, $F_{K,\s}$, $F_{K,\s'}$, $F_{\sigma,\s}$. }
\label{fig_vag_fluxes}
\end{center}
\end{figure}

The VAG discretization is introduced in \cite{Eymard.Herbin.ea:2010} for diffusive problems 
on heterogeneous anisotropic media. Its extension to the hybrid-dimensional 
Darcy flow model is proposed in \cite{BGGM14} based upon the following vector space of degrees of freedom: 
$$
V_\D =\{v_K, v_\s, v_\sigma\in\R, K\in \cells, \s\in \nodes, \sigma\in \faces_\G\}. 
$$
The degrees of freedom are exhibited in Figure \ref{fig_vag_fluxes} 
for a given cell $K$ with one fracture face $\sigma$ in bold. 

\begin{figure}[H]
\begin{center}
  \includegraphics[width=0.4\textwidth]{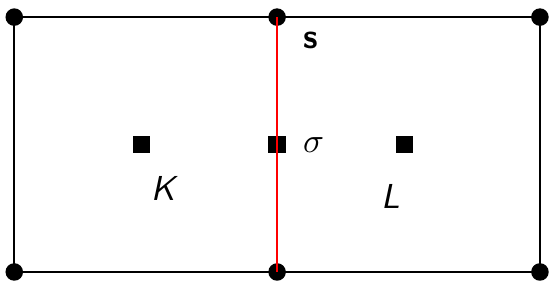}
  \includegraphics[width=0.4\textwidth]{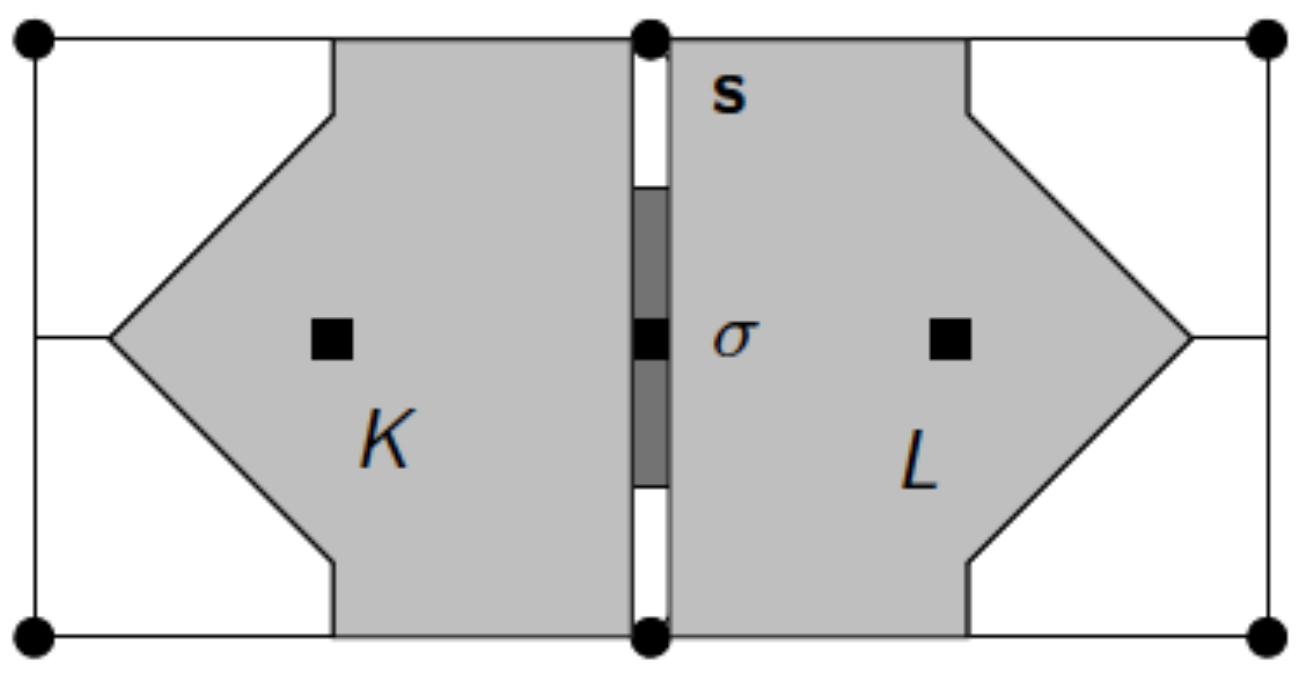}
\caption{{Two cells $K$ and $L$ splitted by 
  one fracture face $\sigma$ in red in the left figure.
  In the right figure, example of control volumes at the two cells $K$ and $L$,
  at the fracture face $\sigma$, and at nodes (the width of the fracture is enlarged in the right figure). 
The control volumes are chosen to avoid to mix the fracture and matrix porous volumes. }}
\label{fig_vag_CV}
\end{center}
\end{figure}

The matrix degrees of freedom are defined by the set of cells $\cells$ and 
by the set of nodes $\nodes\setminus\nodes_\Gamma$ excluding the nodes 
at the matrix fracture interface $\Gamma$. The fracture faces $\faces_\Gamma$ 
and the fracture nodes $\nodes_\Gamma$ are shared between the matrix and the fractures but
the control volumes associated with these degrees of freedom will belong to the fracture network (see Figure \ref{fig_vag_CV}). 
The degrees of freedom at the fracture intersection $\Sigma$ are defined by the set of nodes $\nodes_\Sigma \subset \nodes_\Gamma$ located on $\overline\Sigma$. 
The set of nodes at the Dirichlet boundaries $\overline {\partial \Omega_D}$ 
and $\overline {\partial \Gamma_D}$ is denoted by ${\cal V}_D$. 

{The VAG scheme is a control volume scheme in the sense that it results, 
for each non Dirichlet degree of freedom, in a molar or energy balance equation. }
The matrix diffusion tensor is assumed to be cellwise constant and the tangential 
diffusion tensor in the fracture network is assumed to be facewise constant. 
The two main ingredients are therefore the conservative fluxes and the control volumes. 
The VAG matrix and fracture fluxes are exhibited in Figure \ref{fig_vag_fluxes}. 
For $u_\D\in V_\D$, the matrix 
fluxes $F_{K,\nu}(u_\D)$ connect the cell $K\in\cells$ 
to the degrees of freedom located at the boundary of $K$, namely 
$\nu\in \Xi_K = \nodes_K \cup (\faces_K \cap \faces_\G)$. 
The fracture fluxes $F_{\sigma,\s}(u_\D)$ connect each 
fracture face $\sigma\in \faces_\G$ to its nodes $\s\in \nodes_\sigma$.  
The expression of the matrix (resp. the fracture) fluxes 
is linear and local to the cell (resp. fracture face). 
More precisely, the matrix fluxes are given by 
$$
F_{K,\nu}(u_\D) = \sum_{\nu'\in \Xi_K} T_K^{\nu,\nu'} (u_K - u_{\nu'}), 
$$
with a symmetric positive definite 
transmissibility matrix $T_K = (T_K^{\nu,\nu'})_{(\nu,\nu')\in \Xi_K\times \Xi_K}$ depending 
only on the cell $K$ geometry (including the choices of $\x_K$ and of 
$\x_\sigma, \sigma\in \faces_K$) and on the cell matrix diffusion tensor. 
The fracture fluxes are given by 
$$
F_{\sigma,\s}(u_\D) = \sum_{s\in \nodes_\sigma} T_\sigma^{\s,\s'} (u_\sigma - u_{s'}), 
$$
with a symmetric positive definite 
transmissibility matrix $T_\sigma = (T_\sigma^{\s,\s'})_{(\s,\s')\in \nodes_\sigma\times \nodes_\sigma}$ depending 
only on the fracture face $\sigma$ geometry (including the choice of $\x_\sigma$) 
and on the fracture face width and tangential diffusion tensor. 
Let us refer to \cite{BGGM14} for a more detailed presentation and for the definition 
of $T_K$ and $T_\sigma$. 
 
The construction of the control volumes at each degree of freedom is based on partitions 
of the cells and of the fracture faces. These partitions are respectively denoted, for all $K\in\cells$, by
$$\overline{K} ~ =  ~\overline{\omega}_K ~ \bigcup ~ \left( \bigcup_{\s\in\nodes_K\setminus(\nodes_D\cup \nodes_\Gamma)}\overline{\omega}_{K,\s} 
\right), $$
and, for all $\sigma\in\faces_\G$, by 
$$\overline{\sigma}~=~\overline{\Sigma}_\sigma~\bigcup ~ \left( \bigcup_{\s\in\nodes_\sigma\setminus\nodes_D}\overline{\Sigma}_{\sigma,\s} 
\right). $$
It is important to notice that in the usual case of 
cellwise constant rocktypes in the matrix and facewise constant rocktypes in the fracture network, 
the implementation of the scheme does not require to build explicitly the geometry of these partitions. 
In that case, it is sufficient to define the matrix volume fractions 
$$
\alpha_{K,\s} = {\int_{\omega_{K,\s}} d\x \over \int_K d\x}, 
\s\in\nodes_K\setminus(\nodes_D\cup \nodes_\Gamma), K\in\cells, 
$$
constrained to satisfy $\alpha_{K,\s}\geq 0$,  
and $\sum_{\s\in\nodes_K\setminus(\nodes_D\cup\nodes_\Gamma)}\alpha_{K,\s} \leq 1$, 
as well as the fracture volume fractions 
$$
\alpha_{\sigma,\s} = {\int_{\Sigma_{\sigma,\s}} d_f(\x) d\tau(\x) \over \int_\sigma d_f(\x) d\tau(\x)}, 
\s\in\nodes_\sigma\setminus\nodes_D, \sigma\in\faces_\G, 
$$
constrained to satisfy 
$\alpha_{\sigma,\s}\geq 0$, and $\sum_{\s\in\nodes_\sigma\setminus\nodes_D}\alpha_{\sigma,\s}\leq 1$, where we denote by $d\tau({\bf x})$ the $2$ dimensional Lebesgue measure on $\Gamma$. 
Let us also set 
$$
\phi_K=   (1-\sum_{\s\in \nodes_K\setminus(\nodes_D\cup \nodes_\Gamma)}\alpha_{K,\s})\int_K \phi_m(\x) d\x
 \quad \mbox{ for } K\in \cells, 
$$ 
and 
$$
\phi_\sigma = (1-\sum_{\s\in \nodes_\sigma\setminus\nodes_D}\alpha_{\sigma,\s})
\int_\sigma \phi_f(\x) d_f(\x) d\tau(\x) 
\quad \mbox{ for } \sigma\in \faces_\Gamma, 
$$ 
as well as 
$$
\phi_{\s} =   \sum_{K\in \cells_\s} \alpha_{K,\s}\int_K \phi_m(\x) d\x \quad \mbox{ for } \s \in \nodes\setminus (\nodes_D\cup \nodes_\Gamma), 
$$
and
$$ 
\phi_{\s} = \sum_{\sigma \in \faces_{\Gamma,\s}}\alpha_{\sigma,\s}\int_\sigma \phi_f(\x) d_f(\x) d\tau(\x) 
\quad \mbox{ for } \s \in \nodes_\Gamma\setminus \nodes_D, 
$$
which correspond to the porous volume distributed to the degrees of freedom excluding the  
Dirichlet nodes. 
The rock complementary volumes in each control volume $\nu\in \cells\cup \faces_\Gamma\cup(\nodes\setminus\nodes_D)$ are denoted by $\bar \phi_\nu$. 
 
As shown in \cite{BGGM14}, the flexibility in the choice of the control volumes is a crucial asset, compared with usual CVFE approaches and allows to significantly improve the accuracy of the scheme when the permeability field is highly heterogeneous. As exhibited in Figure \ref{fig_vag_CV}, as opposed to usual 
CVFE approaches, this flexibility allows to define the control volumes 
in the fractures with no contribution from the matrix in order to avoid to enlarge artificially the flow path in the fractures.  

In the following, we will keep the notation $F_{K,s}$, $F_{K,\sigma}$, $F_{\sigma,s}$ for 
the VAG Darcy fluxes defined with the cellwise constant 
matrix permeability $\K_m$ and the facewise constant fracture width $d_f$ and tangential 
permeability $\K_f$. Since the rock properties are fixed,
the VAG Darcy fluxes transmissibility matrices $T_K$ and $T_\sigma$ are computed only once. 

The VAG Fourier fluxes are denoted in the following 
by $G_{K,\s}$, $G_{K,\sigma}$, $G_{\sigma,\s}$. 
They are obtained with the isotropic matrix and fracture thermal conductivities 
averaged in each cell and 
in each fracture face using the previous time step fluid properties.
Hence VAG Fourier fluxes transmissibility matrices need 
to be recomputed at each time step.

\subsection{VAG discretization of the hybrid-dimensional non-isothermal compositional model}

The time integration is based on a fully implicit Euler scheme to avoid severe restrictions on the time 
steps due to the small volumes and high velocities in the fractures. {Note that the thermal conductivities are
discretized as mentioned above using the saturations at the previous time step.} 
A phase based upwind scheme is used for the approximation of the 
mobilities in the Darcy fluxes, that is to say the same scheme that is already used in the 
definition of the transmission conditions \eqref{Hybrid_transmission} of the hybrid-dimensional model. 
At the matrix fracture interfaces, we avoid the mixing of the matrix 
and fracture rocktypes in our choice of the control volumes for $\sigma\in \faces_\Gamma$ 
and $\s\in \nodes_\Gamma$ (see Figure \ref{fig_vag_CV}). 
To avoid too small control volumes at the  
nodes $\s\in \nodes_\Sigma$ located at the fracture intersection, 
the volume is distributed to such a node 
$\s$ from all the fracture faces containing the node $\s$. 
It results that the volumes of the control volumes $\s\in \nodes_\Sigma$ at the fracture intersection is not 
smaller than at any other matrix fracture degrees of freedom. This solves the problems reported in \cite{KDA04} and \cite{SBN12} related to the small 
volumes at the fracture intersections 
and avoid the Star-Delta transformation used in \cite{KDA04} which 
is not valid in the case of multiphase flows. \\

For $N_T\in\N^*$, let us consider the time discretization 
$t^0= 0 < t^1 <\cdots < t^{n-1} < t^n \cdots < t^{N_T} = T$ 
of the time interval $[0,T]$. We denote the time steps by 
$\Delta t^n = t^{n}-t^{n-1}$ for all $n=1,\cdots,N_T$. 

Let be given, for each degree of freedom $\nu \in \cells\cup \faces_\Gamma \cup \nodes$,  
the set of physical unknowns of the Coats' formulation 
$$
X_\nu = \(P_\nu,T_\nu,S_\nu^\alpha, C_\nu^\alpha, \alpha\in Q_\nu, n_{i,\nu}, 
i\in \overline {\cal C}_{Q_\nu}, Q_\nu\).  
$$
We denote by $X_\D$, the full set of unknowns 
$$
X_\D = \{X_\nu, \nu \in \cells\cup \faces_\Gamma \cup \nodes\}. 
$$
We will use the notation $Q_\D = \(Q_\nu, \nu \in \cells\cup \faces_\Gamma \cup \nodes\)$,  
and, for a given $Q_\nu$, we denote by 
$$
X_{Q_\nu} = \(P_\nu,T_\nu,S_\nu^\alpha, C_\nu^\alpha, \alpha\in Q_\nu, n_{i,\nu}, 
i\in \overline {\cal C}_{Q_\nu}\),   
$$
the set of physical unknowns excluding the set of present phases $Q_\nu$. Similarly, for a given $Q_\D$, 
we set  
$$
X_{Q_\D} = \{X_{Q_\nu}, \nu \in \cells\cup \faces_\Gamma \cup \nodes\}. 
$$
We can clearly identify $X_\nu$ and $(X_{Q_\nu},Q_\nu)$ as well as $X_\D$ and $(X_{Q_\D},Q_\D)$.
 
The Darcy fluxes taking into account the gravity term are defined by 
\begin{equation}\label{flux_phase_vag}
\left\{\begin{array}{r@{\,\,}c@{\,\,}ll}
&V^\alpha_{K,\nu}(X_\D) &=& F_{K,\nu}(P_\D) + \rho_{K,\nu}^\alpha F_{K,\nu}({\cal G}_\D), \,\, \nu\in \Xi_K, K\in \cells,\\
&V^\alpha_{\sigma,\s}(X_\D) &=& F_{\sigma,\s}(P_\D) + \rho_{\sigma,s}^\alpha F_{\sigma,\s}({\cal G}_\D), \,\, \s\in \nodes_\sigma, \sigma \in \faces_\Gamma, 
\end{array}\right.
\end{equation}
where ${\cal G}_\D$ denotes the vector $({\bf g}\cdot {\bf x}_\nu)_{\nu \in \cells\cup \faces_\Gamma \cup \nodes}$, 
and the phase mass density is defined by the weighted average 
$$
\rho_{\nu, \nu'}^\alpha = 
{S^{\alpha}_\nu \rho^\alpha(P_\nu,T_\nu,C^\alpha_\nu) + S^{\alpha}_{\nu'}\rho^\alpha(P_{\nu'},T_{\nu'},C^\alpha_{\nu'}) \over S^{\alpha}_{\nu} + S^{\alpha}_{\nu'} }. 
$$
The discretization of the mobilities is obtained using a usual 
phase based upwinding (see e.g. \cite{aziz-settari-79}). For each Darcy flux, 
let us define the phase dependent 
upwind control volume $cv_{\mu,\nu}^\alpha$ such that  
\begin{equation*}
cv_{K,\nu}^\alpha = 
\left\{\begin{array}{r@{\,\,}c@{\,\,}l}
K &\text{ if } & V^\alpha_{K,{\nu}}(X_\D) \geqslant 0 \\
\nu & \text{ if } & V^\alpha_{K,{\nu}}(X_\D) < 0
\end{array}\right. 
\text{ for } K \in \cells, \nu \in \Xi_K,
\end{equation*}
for the matrix fluxes, and such that  
\begin{equation*}
cv_{\sigma,s}^\alpha = 
\left\{\begin{array}{r@{\,\,}c@{\,\,}l}
\sigma &\text{ if } &V^\alpha_{\sigma,s}(X_\D) \geqslant 0 \\
s &\text{ if } &V^\alpha_{\sigma,s}(X_\D) < 0
\end{array}\right.
\text{ for }  \sigma \in \faces_\Gamma, s \in \nodes_\sigma, 
\end{equation*}
for fracture fluxes. Using this upwind discretization, the component molar fluxes are given by 
$$
q_{i,\nu,\nu'}(X_\D) = \sum_{\alpha \in Q_{cv_{\nu,\nu'}^{\alpha}} \cap \mathcal{P}_i} 
m_i^\alpha(X_{cv_{\nu,\nu'}^{\alpha}}) V^\alpha_{\nu,\nu'}(X_\D)
$$
for $i\in {\cal C}$, and the energy fluxes by 
$$
q_{e,\nu,\nu'}(X_\D) = \sum_{\alpha \in Q_{cv_{\nu,\nu'}^{\alpha}} \cap \mathcal{P}_i} 
m_e^\alpha(X_{cv_{\nu,\nu'}^{\alpha}}) V^\alpha_{\nu,\nu'}(X_\D) + G_{\nu,\nu'}(T_\D). 
$$
Next, in each control volume $\nu$, let us denote by   
$$
{\cal A}_{i,\nu}(X_\nu) =  \phi_\nu  n_i(X_\nu), \ i\in {\cal C}
$$
the component molar accumulation,  and by 
$$
{\cal A}_{e,\nu}(X_\nu) =  \phi_\nu  E(X_\nu) + \bar \phi_\nu E_r(P_\nu,T_\nu), 
$$
the energy accumulation. 

We can now state the system of discrete equations at each 
time step $n=1,\cdots,N_T$ which accounts for the component and energy 
conservation equations $i\in {\cal C}\cup \{e\}$ in each cell $K\in \cells$ 
\begin{equation}
\label{VAGconvcell}
R_{K,i} (X_\D^{n}) 
:= \frac{ {\cal A}_i(X_K^{n})-{\cal A}_i(X_K^{n-1})}{\Delta t^n} 
+ \sum_{\s \in \nodes_K}  q_{i,K,\s}(X_\D^n)
+ \sum_{\sigma \in \faces_{\Gamma}\cap\faces_K}  q_{i,K,\sigma}(X_\D^n) = 0,
\end{equation}
in each fracture face $\sigma \in \faces_{\Gamma}$ 
\begin{equation}
\label{VAGconvfrac} 
R_{\sigma,i} (X_\D^{n}) :=
\frac{{\cal A}_i(X_\sigma^{n})-{\cal A}_i(X_\sigma^{n-1})}{\Delta t^n} + 
\sum_{\s \in {\cal V}_\sigma}  q_{i,\sigma,\s}(X_\D^n) 
+ \sum_{K \in {\cal M}_\sigma}  - q_{i,K,\sigma}(X_\D^n) = 0, 
\end{equation}
and at each node $\s\in \nodes\setminus\nodes_D$
\begin{equation}
\label{VAGconvnode} 
R_{s,i} (X_\D^{n}) :=
\frac{{\cal A}_i(X_\s^{n})-{\cal A}_i(X_\s^{n-1})}{\Delta t^n} + 
\sum_{\sigma \in {\cal F}_{\Gamma,\s}}  - q_{i,\sigma,\s} (X_\D^n)
+ \sum_{K \in {\cal M}_\s}  - q_{i,K,\s}(X_\D^n) = 0. 
\end{equation}
It is coupled with the local closure laws 
\begin{equation}
\label{VAG_closure}
\mathbf{0} = L_\nu(X_\nu^n) := 
\left\{\begin{array}{l}
\dsp{ \sum_{\alpha\in Q_\nu^n} S^{\alpha,n}_\nu } - 1,  \\
\dsp{\sum_{i\in {\cal C}^\alpha} C_{i,\nu}^{\alpha,n}} - 1, \,\, \alpha \in Q^n_\nu,\\
f_i^\alpha(P^n_\nu,T_\nu^n,C_\nu^{\alpha,n})C_{i,\nu}^{\alpha,n} - 
f_i^\beta(P^n_\nu,T_\nu^n,C_\nu^{\beta,n})C_{i,\nu}^{\beta,n}, \\  
\quad \quad \quad \quad \quad \alpha\neq\beta,  (\alpha,\beta)\in (Q^n_\nu\cap{\cal P}_i)^2. 
\,\,i\in{\cal C},
\end{array}\right.
\end{equation}
the flash computations $Q^n_\nu =  Q_{flash}(X^n_\nu)$ for 
$\nu \in \cells\cup (\nodes\setminus \nodes_\D) \cup \faces_\Gamma$, 
and the Dirichlet boundary conditions 
$$
X_\s = X_{\s,D}, 
$$
for all $\s\in \nodes_\D$.

\subsection{Newton-Raphson non-linear solver}
\label{sec_newton}

Let us denote by $R_\nu(X_\D)$ the vector $\(R_{\nu,i}, \ i \in \mathcal{C}\cup\{e\}\) $, and 
let us rewrite the conservation equations \eqref{VAGconvcell}, \eqref{VAGconvfrac}, 
\eqref{VAGconvnode} and the closure laws \eqref{VAG_closure} as well 
as the boundary conditions in vector form defining the following 
non-linear system at each time step $n=1,2,...,N_T$
\begin{eqnarray}
\label{NonLinearSystem}
\mathbf{0} = \mathcal{R} (X_\D) := 
\left\{\begin{array}{llllll}
\left(\begin{array}{c}
R_{\s} (X_\D) \\
L_\s (X_s)
\end{array}\right) 
\ \s \in \nodes, \\
\left(\begin{array}{c}
R_{\sigma} (X_\D) \\
L_\sigma (X_\sigma)
\end{array}\right) 
\ \sigma \in \faces_\Gamma, \\
\left(\begin{array}{c}
R_{K} (X_\D) \\
L_K (X_K)
\end{array}\right) 
\ K \in \cells,
\end{array}\right.
\end{eqnarray}
where the superscript $n$ is dropped to simplify the notations 
and where the Dirichlet boundary conditions 
have been included at each Dirichlet node $\s\in {\cal V}_D$ in order to 
obtain a system size independent on the boundary conditions. 

The non-linear system ${\cal R}(X_\D)=0$ coupled to the  flash fixed point 
equations $Q_\nu =  Q_{flash}(X_\nu)$, $\nu\in \cells\cup\faces_\Gamma\cup(\nodes\setminus\nodes_D)$ 
{is solved by an active set Newton-Raphson algorithm widely used in the reservoir simulation community \cite{coats-89-imp} which is detailed below.}
The algorithm is initialized with an initial guess 
$X_{Q_\D}^{(0)}$,$Q_{\D}^{(0)}$ usually given by the previous time step solution and 
computes the initial residual $\mathcal{R}( X_{Q_\D}^{(0)}, Q_{\D}^{(0)})$ 
and its norm $|| \mathcal{R}( X_{Q_\D}^{(0)}, Q_{\D}^{(0)}) ||$ for a given weighted norm $\|.\|$. 

The Newton algorithm iterates on the following steps 
for $r=0,\cdots,$ until convergence of the relative residual 
$$
{|| \mathcal{R}( X_{Q_\D}^{(r)}, Q_{\D}^{(r)}) || \over || 
\mathcal{R}( X_{Q_\D}^{(0)}, Q_{\D}^{(0)}) || } \leq \epsilon_{newton}
$$ 
for a given stopping criteria $\epsilon_{newton}$ or until it reaches a 
maximum number of Newton steps $N_{newton}^{max}$. 

\begin{enumerate}
\item Computation of the Jacobian matrix 
$$
J^{(r)} = {\partial {\cal R} \over \partial X_{{Q_\D}} }
\Big(X_{{Q_\D}}^{(r)},Q^{(r)}_{{\D}} \Big). 
$$
\item Solution of the linear system 
\begin{equation}
\label{solvejac}
J^{(r)} ~d X^{(r)}_{{Q_\D}} = 
-  {\cal R} \Big(X_{{Q_\D}}^{(r)},Q^{(r)}_{{\D}} \Big). 
\end{equation}
\item Update of the unknowns $X_{{Q_\D}}^{(r)}$ 
with a full Newton step $\theta^{(r)}=1$ or a possible relaxation $\theta^{(r)} \in (0,1)$. 
$$
X_{{Q_\D}}^{(r+{1\over 2})} = X_{{Q_\D}}^{(r)} + \theta^{(r)}~d X^{(r)}_{{Q_\D}}.
$$
\item Flash computations to update the sets of present phases $Q_{{\D}}^{(r+1)}$. 
The flash computations also 
provide the molar fractions of the new sets of present phases. 
They are used together with $X_{Q_\D}^{(r+{1\over 2})}$ 
and $Q_{\D}^{(r+1)}$ 
to update the new set of unknowns $X_{Q_\D}^{(r+1)}$. 
\item Computation of the new residual 
$\mathcal{R}( X_{Q_\D}^{(r+1)}, Q_{\D}^{(r+1)})$ and of its norm. 
\end{enumerate}
If the Newton algorithm reaches the predefined maximum number 
of iterations before convergence, 
then we restart this time step with a reduced $\Delta t$.

In view of the non-linear system (\ref{NonLinearSystem}), the size of the linear system for the computation of the Newton step can be considerably reduced without fill-in by 
\begin{itemize}
\item Step 1: elimination of the local closure laws (\ref{VAG_closure}),
\item Step 2: elimination of the cell unknowns. 
\end{itemize}
Step 2 is detailed in Section \ref{sec_parajac}. 
The elimination of the local closure laws (Step 1) is achieved 
for each control volume $\nu\in{\mathcal{M} \cup \mathcal{F}_\Gamma \cup \mathcal{V}} $ by splitting the unknowns $X_{Q_\nu}$ into $\#{\cal C}+1$ primary unknowns 
$X_{Q_\nu}^{pr}$ and $N^{sd}_\nu$ secondary unknowns $X_{Q_\nu}^{sd}$ with 
$$
N^{sd}_\nu =   1 + \# Q_\nu + \sum_{\alpha\in Q_\nu} \#{\cal C}^\alpha 
+\#\overline{\cal C}_{Q_\nu} -\#{\cal C}. 
$$ 
For each control volume $\nu\in{\mathcal{M} \cup \mathcal{F}_\Gamma \cup \mathcal{V}}  $, 
the secondary unknowns are  chosen in such a way that the square matrix 
$$
{\partial L \over \partial X^{sd}_{Q_\nu} } \left(X^{pr}_{Q_\nu},X_{Q_\nu}^{sd}, Q_\nu\right)
\in \mathbb{R}^{N^{sd}_\nu \times N^{sd}_\nu},  
$$
is non-singular. This choice can be done algebraically in the general case, 
or defined once and for all for each set of present phases $Q_\nu$ 
for specific physical cases. Here we remark that the unknowns $(n_{i,\nu})_{i\in\overline{\cal C}_{Q_\nu}}$ are not involved in the closure laws \eqref{VAG_closure} and hence are 
always chosen as primary unknowns.

The ill conditioned linear system obtained from \eqref{solvejac} after the two elimination steps is solved using an iterative solver such as GMRES or BiCGStab 
combined with a preconditioner adapted to the elliptic 
or parabolic nature of the pressure unknown and to the coupling with the remaining hyperbolic or parabolic unknowns. 
One of the most efficient preconditioners for such systems is the so-called CPR-AMG preconditioner 
introduced in \cite{LVW:2001} and \cite{SMW:2003}. 
It combines multiplicatively an algebraic multigrid preconditioner (AMG) for a pressure block of the linear system \cite{HY:2001} 
with a more local preconditioner for the full system, 
such as an incomplete LU factorization. 
The choice of the pressure block is important for the efficiency of the CPR-AMG preconditioner. In the following 
experiments we simply define the pressure equation in each control volume by the sum of the molar conservation equations 
in the control volume. Let us refer to \cite{LVW:2001}, \cite{SMW:2003}, and \cite{ABMQ:2006} for a discussion 
on other possible choices. 
Let us denote by $Jx=b$ the linear system where $J$ is the Jacobian matrix 
and $b$ the right hand side taking into account 
the elimination steps and the linear combinations of the lines for the pressure block. 
The CPR-AMG preconditioner $\mathcal{P}_{\mathrm{cpr-amg}}$ is defined for any vector $b$ 
by $\mathcal{P}_{\mathrm{cpr-amg}} b = v$ with 
\begin{gather}
v^{1/2} = R_{P}^{\top} \mathcal{P}_{\mathrm{amg}(J_P)} R_{P} b,  \nonumber \\
v = v^{1/2} + \mathcal{P}_{\mathrm{ilu0}} (b - J v^{1/2}), 
\end{gather}
where  $\mathcal{P}_{\mathrm{amg}(J_P)}$ is the AMG preconditioner with $ J_P = R_{P} J R_{P}^{\top}$, $\mathcal{P}_{\mathrm{ilu0}}$ is the ILU(0) preconditioner applied on the Jacobian $J$, $R_{P}$ is the restriction matrix to the pressure unknowns and $R_{P}^{\top}$ is the transpose of $R_{P}$.


\section{Parallel implementation}
\label{sec_pl}

Parallel implementation is achieved using the Message Passing Interface (MPI). Let us denote by $N_p$ the number of MPI processes.

\subsection{Mesh decomposition}

The set of cells $\cells$ is partitioned into $N_p$ subsets $\cells^p, p=1,...,N_p$ using the library METIS \cite{citemetis}.
{In the current implementation, this partitioning is only based on the cell connectivity graph and does not take into account the fracture faces.
This will be investigated in the near future and the potential gain is discussed in the numerical section. } 
The partitioning of the set of nodes $\mathcal{V}$ and of the set of fracture faces $\mathcal{F}_{\Gamma}$ is defined as follows:
assuming we have defined a global index of the cells $K\in \cells$ 
let us denote 
by $K(\s), \s\in\nodes$ (resp. $K(\sigma)$, $\sigma\in \faces_\Gamma$) the cell with the smallest global index among those of $\cells_\s$ 
(resp. $\cells_\sigma$). Then we set 
$$
\nodes^p = \{\s\in\nodes \,|\, K(\s)\in \cells^p\}, 
$$
and 
$$
\faces_\Gamma^p = \{\sigma\in\faces_\Gamma \,|\, K(\sigma)\in \cells^p\}. 
$$

\begin{figure}[H]
\centering
\includegraphics[width=0.4\textwidth]{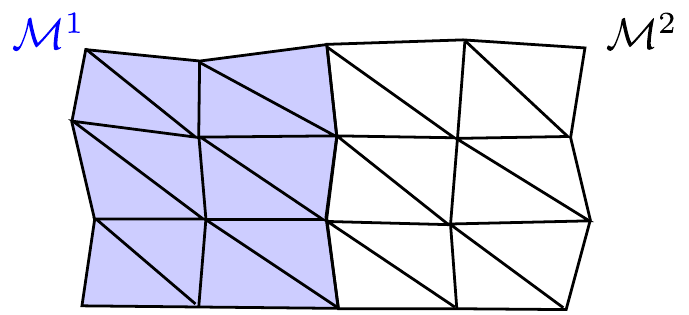}
\hspace{0.05\textwidth}
\includegraphics[width=0.355\textwidth]{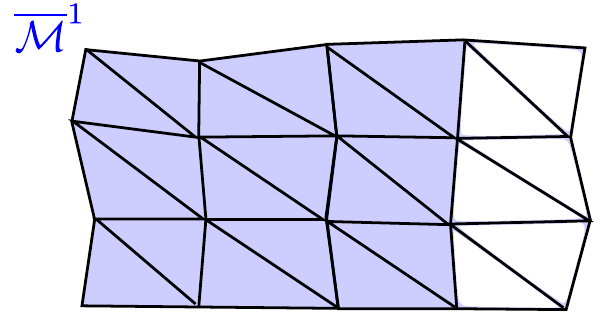} \\
\hspace{0.05\textwidth}
\includegraphics[width=0.4\textwidth]{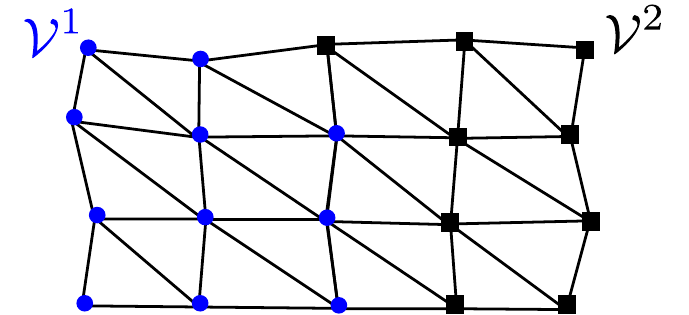}
\hspace{0.04\textwidth}
\includegraphics[width=0.41\textwidth]{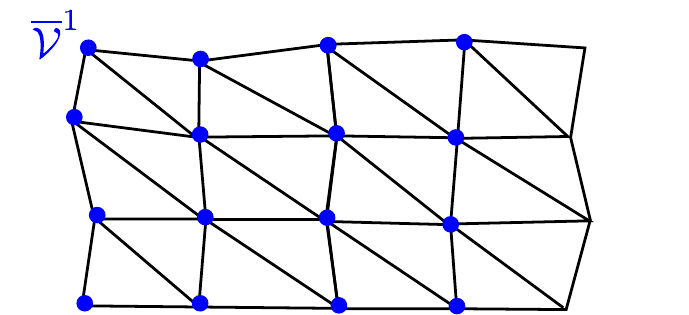} 
\caption{Example of mesh decomposition.}
\label{fig_cellghost} 
\end{figure}

The overlapping decomposition of $\cells$ into the sets
$$
\overline{\cells}^p, \,\, p=1,...,N_p, 
$$
is chosen in such a way that any compact finite volume scheme such as the VAG scheme can be assembled locally on each process. 
Hence, as exhibited in Figure~\ref{fig_cellghost}, $\overline{\cells}^p$ is defined as the set of cells sharing a node with a cell of $\cells^p$. 
The overlapping decompositions of the set of nodes and of 
the set of fracture faces follow from this definition:
$$
\overline\nodes^p = \bigcup_{K\in \overline\cells^p} \nodes_K, \,\, p=1,\cdots,N_p, 
$$
and 
$$
\overline\faces_\Gamma^p = \bigcup_{K\in \overline\cells^p} \faces_K\cap\faces_\Gamma, \,\, p=1,\cdots,N_p.  
$$

The partitioning of the mesh is performed by the master process (process 1), and then, each local mesh is distributed to its process.
Therefore, each MPI process contains the local mesh 
($\overline{\cells}^p$, $\overline\nodes^p$,  $\overline\faces_\Gamma^p$), $p=1,2,...,N_p$ which
is splitted into two parts: 
\begin{equation*}
\begin{split}
& \text{own mesh: }  (\cells^p, \nodes^p, \faces_\Gamma^p), \\
& \text{ghost mesh: }  (\overline{\cells}^p \backslash \cells^p, \overline\nodes^p \backslash \nodes^p, \overline\faces_\Gamma^p \backslash \faces_\Gamma^p).
\end{split}
\end{equation*}

We now turn to the parallel implementation of the Jacobian system.

\subsection{Parallelization of the Jacobian system}
\label{sec_parajac}

On each process $p=1,...,N_p$, the local Jacobian system is defined 
by the set of unknowns $X_\nu$, $\nu\in \overline\nodes^p\cup \overline\faces_\Gamma^p\cup\overline \cells^p$, 
the closure equations on control volume $\nu\in \overline\nodes^p\cup \overline\faces_\Gamma^p\cup\overline \cells^p$ 
and the conservation equations of all own nodes $\s \in \mathcal{V}^p$, all own fracture faces $\sigma \in \mathcal{F}_\Gamma^p$ 
and all own and ghost cells $ k\in \overline{\cells}^p$. \\

The local Jacobian system is firstly reduced by eliminating the local closure laws on each control volume $\nu\in \overline\nodes^p\cup \overline\faces_\Gamma^p\cup\overline \cells^p$ using the procedure presented in Section \ref{sec_newton}. 
The local reduced Jacobian system can be written as the following rectangular linear system 
\begin{equation*}
\begin{pmatrix}
J_{s s}^p & J_{s f}^p& J_{s c}^p \\
J_{f s}^p & J_{f f}^p& J_{f c}^p \\
J_{c s}^p & J_{c f}^p& J_{c c}^p
\end{pmatrix}
\begin{pmatrix}
\overline U_s^p \\
\overline U_f^p \\
\overline U_c^p
\end{pmatrix}
=
\begin{pmatrix}
b_s^p \\
b_f^p \\
b_c^p
\end{pmatrix}
\end{equation*} 
where $\overline U_s^p \in \mathbb{R}^{ \# \overline{\mathcal{V}}^p \times (\# \mathcal{C}+1)}$, 
$\overline U_f^p \in \mathbb{R}^{\# \overline{\mathcal{F}}_\Gamma^p \times (\# \mathcal{C}+1)}$ 
and $\overline U_c^p \in \mathbb{R}^{\# \overline{\cells}^p \times (\# \mathcal{C}+1)}$ denote the vector of own and ghost primary unknowns $X^{pr}_{Q_\nu}$
at the nodes $\nu\in {\overline \nodes}^p$, at the 
fracture faces $\nu\in {\overline \faces}_\Gamma^p$ and at the 
cells $\nu\in {\overline \cells}^p$ respectively on the process $p$. 
The above matrices have the following sizes 
\begin{align*}
& J_{ss}^p    \in \mathbb{R}^{\big( \#\mathcal{V}^p \times (\# \mathcal{C}+1) \big) \times  \big( \# \overline{\mathcal{V}}^p \times (\# \mathcal{C}+1) \big)}, \\
& J_{ff}^p    \in \mathbb{R}^{\big( \# \mathcal{F}_{\Gamma}^p \times (\# \mathcal{C}+1) \big) \times  \big( \# \overline{\mathcal{F}}_{\Gamma}^p \times (\# \mathcal{C}+1) \big)}, \\
& J_{cc}^p    \in \mathbb{R}^{\big( \# \overline{\cells}^p \times (\# \mathcal{C}+1) \big) \times  \big( \# \overline{\cells}^p \times (\# \mathcal{C}+1) \big)},
\end{align*}
and $b_s^p \in \mathbb{R}^{ \# {\mathcal{V}}^p \times (\# \mathcal{C}+1)}$, 
$b_f^p \in \mathbb{R}^{\# {\mathcal{F}}_\Gamma^p \times (\# \mathcal{C}+1)}$ 
and $b_c^p \in \mathbb{R}^{\# \overline{\cells}^p \times (\# \mathcal{C}+1)}$
denote the corresponding right hand side vectors. The matrix $J_{cc}^p$ is a non singular diagonal matrix and the cell unknowns can be 
easily eliminated without fill-in leading to the following Schur complement system 
\begin{equation}
\label{eqschur}
J^p 
\begin{pmatrix}
\overline U_s^p \\
\overline U_f^p
\end{pmatrix}
=
b^p,
\end{equation}
with
\begin{equation*}
\begin{split}
& J^p:=
\begin{pmatrix}
J_{ss}^p & J_{sf}^p \\
J_{fs}^p & J_{ff}^p
\end{pmatrix}
-
\begin{pmatrix}
J_{sc}^p \\
J_{fc}^p
\end{pmatrix}
(J_{cc}^p)^{-1}
\begin{pmatrix}
J_{cs}^p & J_{cf}^p
\end{pmatrix}
,\\
& b^p :=
\begin{pmatrix}
b_s^p \\
b_f^p
\end{pmatrix}
-
\begin{pmatrix}
J_{sc}^p \\
J_{fc}^p
\end{pmatrix}
(J_{cc}^p)^{-1} b_c^p,
\end{split}
\end{equation*}
and
\begin{equation}
\label{eqschurcell}
\overline U_c^p = (J_{cc}^p)^{-1} (b_c^p - J_{cs}^p \overline U_s^p - J_{cf}^p \overline U_f^p).
\end{equation}
The linear system \eqref{eqschur} is built locally on each process $p$ and
transfered to the parallel linear solver library PETSc \cite{citepetsc}.
The parallel matrix and the parallel vector in PETSc are stored in a distributed manner, i.e. each process stores its own rows. We construct the following parallel global linear system 
\begin{equation}
\label{parallellin}
J U= b,
\end{equation}
with
\begin{equation*}
J:=
\begin{pmatrix}
J^1 R^1 \\
J^2 R^2 \\
\vdots \\
J^{N_p} R^{N_p}
\end{pmatrix}
\begin{array}{l}
\big\} \text{ process 1}  \\
\big\} \text{ process 2}  \\
\quad \quad \vdots \\
\big\} \text{ process $N_p$}
\end{array},
\end{equation*}
and
\begin{equation*}
U:=
\begin{pmatrix}
U_s^1 \\
U_f^1 \\
U_s^2 \\
U_f^2 \\
\vdots \\
U_s^{N_p} \\
U_f^{N_p}
\end{pmatrix}
\begin{array}{l}
\bigg\} \text{ process 1} \\
\bigg\} \text{ process 2} \\
\quad \quad \vdots \\
\bigg\} \text{ process $N_p$}
\end{array}, \
b:=
\begin{pmatrix}
b^1 \\
b^2 \\
\vdots \\
b^{N_p}
\end{pmatrix}
\begin{array}{l}
\big\} \text{ process 1 }  \\
\big\} \text{ process 2 }  \\
\quad \quad \vdots \\
\big\} \text{ process $N_p$}
\end{array}
\end{equation*}
where $R^p,p=1,2,...,N_p$ is a restriction matrix satisfying 
$$R^p U = \begin{pmatrix}
\overline U_s^p \\
\overline U_f^p
\end{pmatrix}.$$
The matrix $J^p R^p$, the vector $\begin{pmatrix}
U_s^p \\
U_f^p 
\end{pmatrix}$
and the vector $b^p$ are stored in process $p$.

The linear system \eqref{parallellin} is solved using 
the GMRES algorithm preconditioned by CPR-AMG preconditioner
as discussed in the previous section. 
The solution of the linear system provides on each process $p$ the solution vector $\begin{pmatrix}
U_s^p \\
U_f^p 
\end{pmatrix}$
of own node and fracture-face unknowns. 
Then, the ghost node unknowns $U^p_\nu$, $\nu \in (\overline\nodes^p \backslash \nodes^p)$ 
and the ghost fracture face unknowns $U_\nu^p$, $\nu\in  (\overline\faces_\Gamma^p \backslash \faces_\Gamma^p)$ are recovered 
by a synchronization step with MPI communications. This synchronization is efficiently 
implemented using a PETSc  matrix vector product 
\begin{equation}
\label{commS}
\overline U = S U
\end{equation}
where
$$
\overline U:=
\begin{pmatrix}
\overline U_s^1 \\
\overline U_f^1 \\
\overline U_s^2 \\
\overline U_f^2 \\
\vdots
\end{pmatrix}
$$
is the vector of own and ghost node and fracture-face unknowns on all processes. The matrix 
$S$, containing only $0$ and $1$ entries, is assembled once and for all at the beginning of the simulation. 


Finally, thanks to \eqref{eqschurcell}, the vector of own and ghost cell unknowns $\overline U_c^p$ is 
computed locally on each process $p$.

\section{Numerical results}
\label{sec_num}

The numerical tests are {all} implemented in the framework of the code ComPASS 
on the cluster ``cicada" hosted by University Nice Sophia-Antipolis 
consisting of 72 nodes (16 cores/node, Intel Sandy Bridge E5-2670, 64GB/node). 
We always fix 1 core per process and 16 processes per node. 
The communications are handled by OpenMPI 1.8.2 (GCC 4.9). 

Five test cases are considered in the following subsections. They include 
a two-phase immiscible isothermal Darcy flow model, a two-phase isothermal Black Oil model and a non-isothermal liquid gas flow model. 
Different types of meshes namely hexahedral, tetrahedral, prismatic and Cartesian meshes are used in these simulations.

The settings of the nonlinear Newton and linear GMRES solvers 
are defined by their maximum number of iterations denoted by $N_{newton}^{max}$ and $N_{gmres}^{max}$ and 
by their stopping criteria on the relative residuals denoted by $\epsilon_{newton}$ and $\epsilon_{gmres}$. 

The time stepping is defined by an  
initial time step $\Delta t^{(0)}$ and by a maximum time step $\Delta t^{(k)}$ on each time interval 
$[t^{(k)},t^{(k+1)})$, $k=0,\cdots,k_f-1$ with $t^{(0)} = 0$ and $t^{(k_f)} = t_f$, where $t_f$ is the final 
simulation time. The successive time steps are computed using the following rules. 
If the Newton algorithm 
reaches convergence in less than $N_{newton}^{max}$ iterations at time step $n$ with $t^n \in [t^{(k)},t^{(k+1)})$,  
then the next time step $\Delta t^{n+1}$ is set to 
\begin{equation}
\label{choosedt}
\Delta t^{n+1} = \mathrm{min}(c \Delta t^n, \Delta t^{(k)}),\ c =1.2.
\end{equation} 
If the Newton algorithm does not converge in $N_{newton}^{max}$ iterations or if 
the linear solver does not reach convergence in $N_{gmres}^{max}$  iterations, 
then the time step is chopped by a factor two and restarted.

In all the following numerical experiments, 
the relative permeabilities are given by the Corey laws $k_r^{\alpha}(S)=(S^{\alpha})^2$ 
for both phases $\alpha \in \mathcal{P}$ and both in the matrix domain and in the fracture network.

%

\subsection{Two-phase immiscible isothermal flow}
\label{Sec_num_immiscible}

In this subsection, we consider an immiscible isothermal two-phase  Darcy flow with $\mathcal{P}=\{ \mathrm{water, oil}  \}$ the set of phases 
and $\mathcal{C}=\{ H_2O, HC \}$ the set of components. The model prescribes the mass conservation and we set  
$\rho^{water} = 1000 \text{ kg/m$^3$}$ and $\rho^{oil}= 700 \text{ kg/m$^3$}$. The phase viscosities are  
set to $\mu^{water} = 10^{-3}$ $\text{Pa}\cdot \text{s}$ and $\mu^{oil} = 5.0\times 10^{-3}$ $\text{Pa}\cdot \text{s}$. 

\begin{figure}[!htbp]
  \centering
  \includegraphics[width=0.35\textwidth]{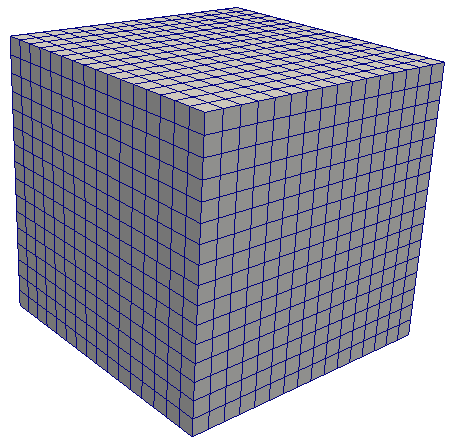}
  \hspace{0.1\textwidth}
  \includegraphics[width=0.35\textwidth]{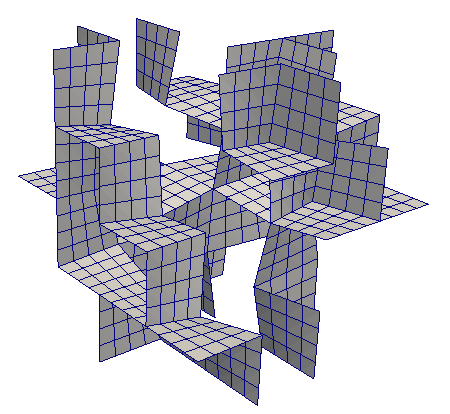}
  \caption{Hexahedral mesh of the matrix domain (left) 
    conforming to the fracture network (right) obtained with $n_x=16$.}
  \label{meshcpgfrac}
\end{figure}

The reservoir domain is defined by $\Omega=(0,100)^3$ in meter.
{We consider a topologically Cartesian mesh of size $n_x\times n_x\times n_x$ of the domain $\Omega$. The mesh is  
exhibited in Figure \ref{meshcpgfrac} for $n_x = 16$.}
The mesh is exponentially refined at the interface between the matrix domain and the fracture network 
as shown in Figure \ref{meshcpgfrac}. 
The width of the fractures is fixed to $d_f=0.01$ meter. The permeabilities are isotropic and set to $\Lambda_m=10^{-15}$ m$^2$ 
in the matrix domain and to $\Lambda_f=10^{-11}$ m$^2$ in the fracture network. The porosities in the matrix domain and in the fractures are $\phi_m=0.1$ and $\phi_f=0.5$ respectively.

The reservoir is initially saturated with water and oil is injected at the bottom boundaries of the matrix domain 
and of the fracture network. The oil phase rises by gravity in the matrix and in the fracture network. 
The lateral boundaries are considered impervious.
The initial pressure is hydrostatic with $P=2$ MPa at the bottom boundaries and $P=1$ MPa at the top boundaries.

The linear and nonlinear solver parameters are fixed to $N_{newton}^{max}=35$, 
$N_{gmres}^{max}=150$, $\epsilon_{gmres}= 10^{-4}$, $\epsilon_{newton}= 10^{-5}$, 
and the time stepping parameters are fixed to $t_f=10000$ days, 
$\Delta t^{(0)} = 5 \text{ days}$, $\Delta t^{(1)} = 5 \text{ days}$, $\Delta t^{(2)} = 15 \text{ days}$, 
$t^{(1)} = 100 \text{ days}$, $k_f = 2$. 

\begin{figure}[!htbp]
  \centering
  \addtocounter{subfigure}{2} 
  \begin{subfigure}{0.35\textwidth}
    \centering
    \includegraphics[width=\textwidth]{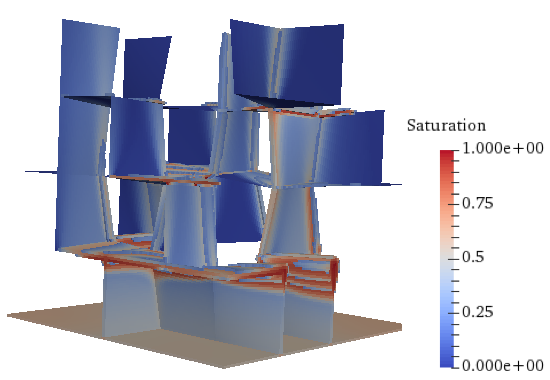}
    \caption{$t=2500$}
  \end{subfigure}
  \hspace{0.05\textwidth}
  \begin{subfigure}{0.35\textwidth}
    \centering
    \includegraphics[width=\textwidth]{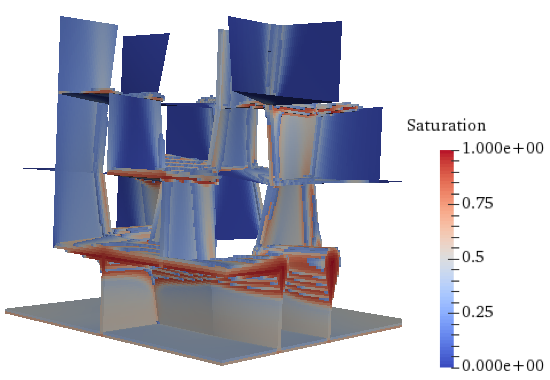}
    \caption{$t=5000$}
  \end{subfigure}
  \begin{subfigure}{0.35\textwidth}
    \centering
    \includegraphics[width=\textwidth]{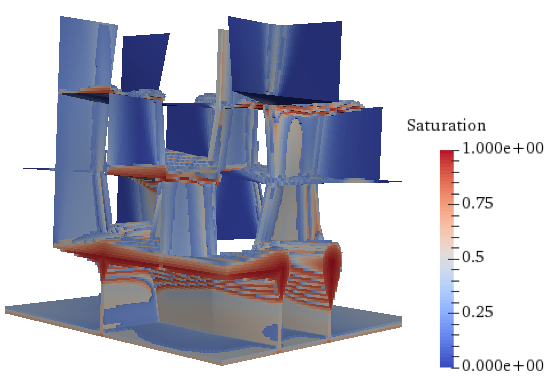}
    \caption{$t=7500$}
  \end{subfigure}
  \hspace{0.05\textwidth}
  \begin{subfigure}{0.35\textwidth}
    \centering
    \includegraphics[width=\textwidth]{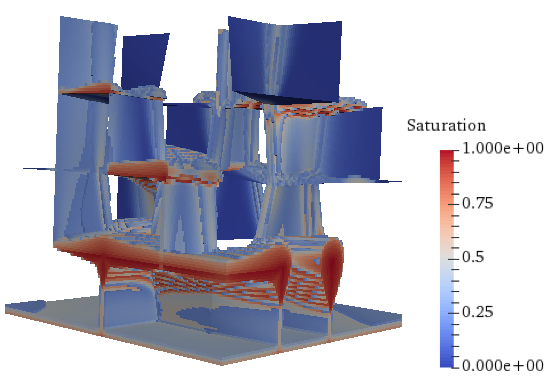}
    \caption{$t=10000$}
  \end{subfigure}
  \caption{Saturation of oil in the fractures and in the matrix domain at different times (in days) for the hexahedral mesh with $n_x=128$. 
A threshold of 0.2 is used for the saturation in the matrix domain.}
  \label{simucpg}
\end{figure}

Figure \ref{simucpg} exhibits the oil saturation obtained with the mesh size $n_x=128$
at times $t=2500, 5000, 7500$ days and at the final time $t_f=10000$ days. 

\begin{table}[!htbp]
\centering
\caption{{Maximum/mean number of own cells, own cells+nodes+fracture faces and own nodes+fracture faces by process for the hexahedral mesh with $n_x=128$ and $N_p=64,128$.} }
\label{partitioning_cpg}
\begin{tabular}{c|c|c}
\hline
$N_p$ & 64 & 128  \\
\hline
own cells & 32768/32768 & 16385/16384    \\
own cells+nodes+fracture faces & 68985/67126 &  34831/33563   \\
own nodes+fracture faces & 36217/34358 & 18447/17179    \\
\hline
\end{tabular}
\end{table}

Table \ref{Nitercpg} clearly shows that both the total numbers of Newton iterations and of linear solver iterations 
are almost independent on the number of MPI processes. The Newton solver requires an average of $2.6$ iterations per time step  
and the GMRES linear solver converges in an average of 40 iterations. These results are very good 
given the mesh size combined with the large constrast of permeabilities and of space and time scales between the fracture network 
and the matrix.  

\begin{table}[!htbp]
\centering
\caption{Number of time steps ($N_{timestep}$), total number of Newton iterations ($N_{newton}$) and total number of linear solver iterations ($N_{gmres}$) vs. number of MPI processes for the two-phase immiscible isothermal flow test case with the hexahedral mesh obtained for $n_x=128$.}
\label{Nitercpg}
\begin{tabular}{c|c|c|c|c}
\hline
$N_p$ & 16 & 32 & 64 & 128 \\
\hline
$N_{timestep}$ & 683 & 683   & 683   & 683  \\
$N_{newton}$ & 1743  & 1742  & 1745  & 1741 \\
$N_{gmres}$ & 68779 & 69015 & 68927 & 69070 \\
$N_{newton}/N_{timestep}$ & 2.6 & 2.6 & 2.6 & 2.5 \\
$N_{gmres}/N_{newton}$ & 39.5 & 39.6 & 39.5 & 39.7 \\
\hline
\end{tabular}
\end{table}

\begin{figure}[!htbp]
  \centering
  \includegraphics[width=0.5\textwidth]{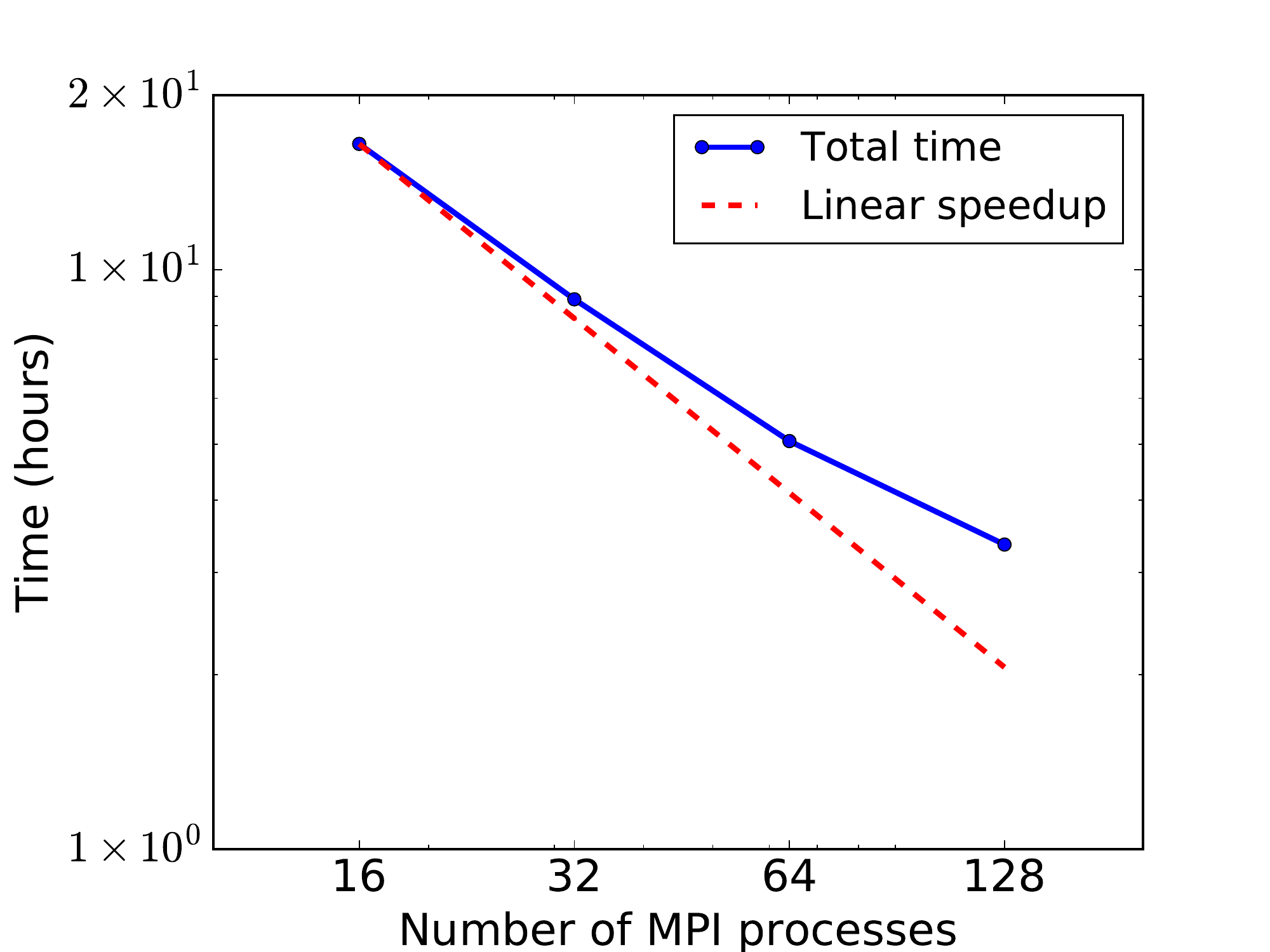}
  \caption{Total computation time vs. number of MPI processes for the two-phase immiscible isothermal flow test case with the
    hexahedral mesh obtained for $n_x=128$.}
  \label{timecpg}
\end{figure}

Figure \ref{timecpg} presents the total computation times in hours for different number of MPI processes $N_p=16,32,64,128$. 
The scalability behaves as expected for fully implicit time integration and AMG type preconditioners. 
It is well known that the AMG preconditioner requires a sufficient number of unknowns per MPI process, 
say $100000$ as classical order of magnitude,  to achieve a linear strong scaling. 
For this mesh size, leading to roughly $2\times 10^6$ unknowns for the pressure block, 
the scalability is still not far from linear on up to $64$ processes and then degrades more rapidly for $N_p = 128$.  {Table \ref{partitioning_cpg} shows that the partitioning could be improved by using a weighted graph taking into account the fracture faces. Nevertheless, for this test case, the potential gain seems rather small compared with the loss of parallel efficiency exhibited in Figure \ref{timecpg} which is mainly due to the communication overhead.} 

\subsection{Black Oil model}

\subsubsection{Oil migration}
This test case considers a Black Oil model with two components $\mathcal{C}=\{ H_2O, HC \}$ 
and two phases $\mathcal{P}=\{ \mathrm{water, oil}  \}$. The $HC$ component can dissolve 
in the water phase defined as 
a mixture of $HC$ and $H_2O$ while the oil phase contains only the $HC$ component. 
The viscosities of the water and oil phases are the same as in the previous test case. 
The mass densities are defined by 
\begin{equation*}
\rho^{water} = 990 \times (1 + C_{HC}^{water}) \text{ kg/m$^{3}$}, \
\rho^{oil} = 700 \text{ kg/m$^{3}$}.
\end{equation*}
The fugacity coefficients $f_{HC}^\alpha(P,T,C^\alpha), \alpha \in \mathcal{P}$ are defined by 
\begin{equation*}
\begin{split}
& f_{HC}^{water} = 1, \\
& f_{HC}^{oil} = \frac{P-P_2}{P_1-P_2} \bar c_1  +  \frac{P-P_1}{P_2-P_1} \bar c_2,
\end{split}
\end{equation*}
with $P_1=1$ MPa, $P_2=2$ MPa, and $\bar c_1 = 5 \times 10^{-3}$, $\bar c_2 = 10^{-2}$. 

\begin{figure}[!htbp]
  \centering
  \includegraphics[width=0.4\textwidth]{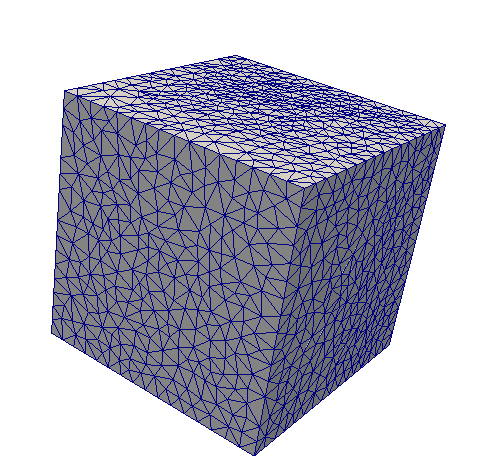}
  \hspace{0.1\textwidth}
  \includegraphics[width=0.4\textwidth]{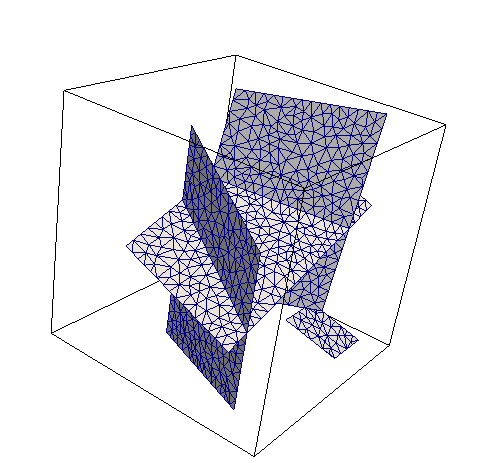}
  \caption{Example of tetrahedral mesh of the matrix domain (left) 
    conforming to the fracture network (right).}
  \label{meshtet}
\end{figure}

The reservoir is the cubic domain $\Omega=(0,100)^3$ in meter and the width of the fractures is fixed to 
$d_f=0.01$ meter. We consider a tetrahedral mesh conforming to the fracture network as 
exhibited in Figure~\ref{meshtet} for a coarse mesh. 
The mesh used in this subsection contains about $6.2\times 10^6$ cells, $9.7\times 10^5$ nodes and $7.1\times 10^4$ fracture faces.
The permeabilities are isotropic and fixed to $\Lambda_m=10^{-15} $ m$^2$ in the matrix domain and to 
$\Lambda_f=10^{-11} $ m$^2$ in the fracture network. The porosities in the matrix domain and in the fractures are $\phi_m=0.1$ and $\phi_f=0.5$ respectively.

As in the previous test case, the reservoir is initially saturated with pure water and oil is injected at 
the bottom boundaries of the matrix domain and of the fracture network. The initial pressure is hydrostatic with $P=2$ MPa at the bottom boundaries and $P=1$ MPa at the top boundaries. 

The linear and nonlinear solver parameters are fixed to $N_{newton}^{max}=35$, 
$N_{gmres}^{max}=200$, $\epsilon_{gmres}= 10^{-4}$, $\epsilon_{newton}= 10^{-5}$, 
and the time stepping parameters are fixed to $t_f=10000$ days, 
$\Delta t^{(0)} = 0.5 \text{ days}$, $\Delta t^{(1)} = 2 \text{ days}$, 
$\Delta t^{(2)} = 50 \text{ days}$, $\Delta t^{(3)} = 100 \text{ days}$, 
$t^{(1)} = 180 \text{ days}$, $t^{(2)} = 2000 \text{ days}$, $k_f = 3$.

\begin{figure}[!htbp]
  \centering
  \addtocounter{subfigure}{2} 
  \begin{subfigure}{0.35\textwidth}
    \centering
    \includegraphics[width=\textwidth]{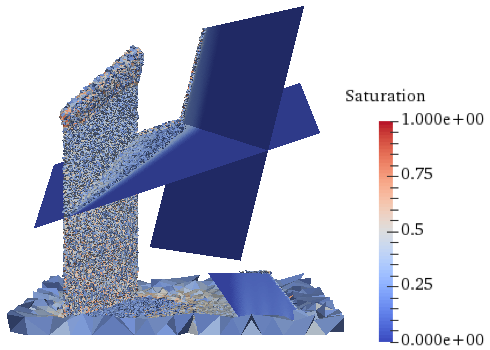}
    \caption{$t=2500$}
  \end{subfigure}
  \hspace{0.05\textwidth}
  \begin{subfigure}{0.35\textwidth}
    \centering
    \includegraphics[width=\textwidth]{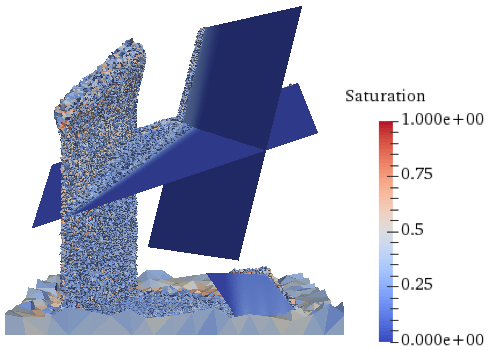}
    \caption{$t=5000$}
  \end{subfigure}
  \begin{subfigure}{0.35\textwidth}
    \centering
    \includegraphics[width=\textwidth]{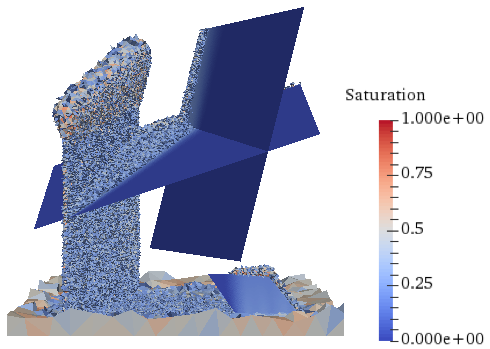}
    \caption{$t=7500$}
  \end{subfigure}
  \hspace{0.05\textwidth}
  \begin{subfigure}{0.35\textwidth}
    \centering
    \includegraphics[width=\textwidth]{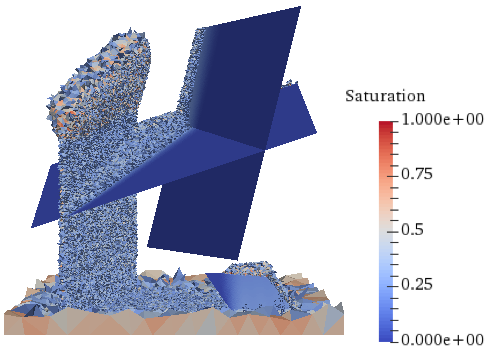}
    \caption{$t=10000$}
  \end{subfigure}
  \caption{Oil saturation in the fractures and in the matrix domain 
    at different times (in days) for the tetrahedral mesh. 
    A threshold of 0.2 is used for the saturation in the matrix domain.}
  \label{simutet_Sat}
\end{figure}

\begin{figure}[!htbp]
  \centering
  \addtocounter{subfigure}{2} 
  \begin{subfigure}{0.35\textwidth}
    \centering
    \includegraphics[width=\textwidth]{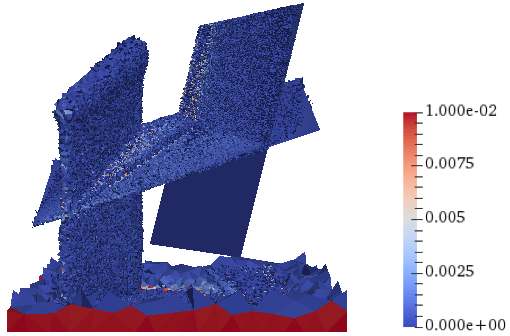}
    \caption{$t=2500$}
  \end{subfigure}
  \hspace{0.05\textwidth}
  \begin{subfigure}{0.35\textwidth}
    \centering
    \includegraphics[width=\textwidth]{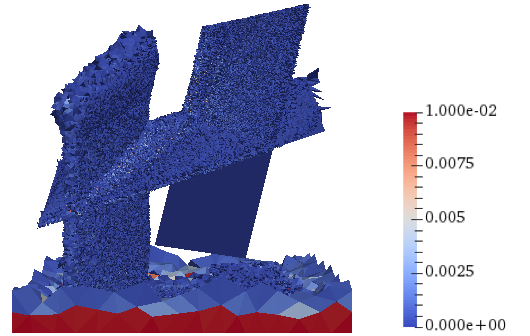}
    \caption{$t=5000$}
  \end{subfigure}
  \begin{subfigure}{0.35\textwidth}
    \centering
    \includegraphics[width=\textwidth]{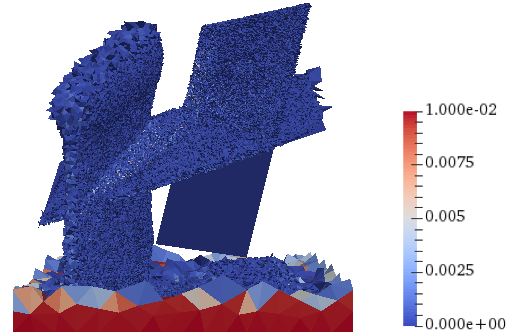}
    \caption{$t=7500$}
  \end{subfigure}
  \hspace{0.05\textwidth}
  \begin{subfigure}{0.35\textwidth}
    \centering
    \includegraphics[width=\textwidth]{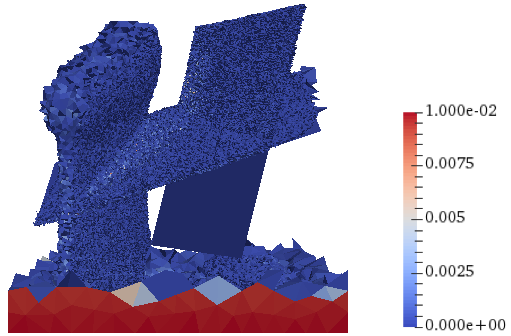}
    \caption{$t=10000$}
  \end{subfigure}
  \caption{Molar fraction of the $HC$ component in the water phase in the fractures and in the matrix domain at different times (in days) 
for the tetrahedral mesh. A threshold of 0.001 is used for the saturation in the matrix domain.}
  \label{simutet_HC}
\end{figure}

Figure \ref{simutet_Sat} and Figure \ref{simutet_HC} present the oil saturation and the molar fraction of the $HC$ component 
in the water phase both in the fractures and in the matrix domain at times $t=2500, 5000, 7500, 10000$ days.

\begin{table}[!htbp]
\centering
\caption{{Maximum/mean number of own cells, own cells+nodes+fracture faces and own nodes+fracture faces by process for the tetrahedral mesh with $6.2\times 10^6$ cells and $N_p=64,128$.} }
\label{partitioning_tet}
\begin{tabular}{c|c|c}
\hline
$N_p$ & 64 & 128  \\
\hline
own cells & 96518/96517 &  48260/48258 \\
own cells+nodes+fracture faces &  114898/112762 &  58063/56381\\
own nodes+fracture faces & 18381/16246  & 9804/8123 \\
\hline
\end{tabular}
\end{table}

As in the previous test case, table \ref{Nitertet} exhibits 
that both the total numbers of Newton iterations and of linear solver iterations are 
almost independent on the number of MPI processes. The average number of Newton iteration is $8.8$ per time step. 
This is a significant increase compared with the previous test case which is due to the 
phase appearance and disappearance in the Black oil model combined with large contrasts of permeabilities 
and space and time scales between the matrix and the fractures. On the other hand, the average number of linear solver iterations 
is roughly $30$ per Newton step which is even better than in the previous test case. 
\begin{table}[!htbp]
\centering
\caption{Number of time steps ($N_{timestep}$), total number of Newton iterations ($N_{newton}$) and total number of linear solver iterations ($N_{gmres}$) vs. number of MPI processes for the black oil model test case with $6.2\times 10^6$ tetrahedral cells.}
\label{Nitertet}
\begin{tabular}{c|c|c|c|c|c}
\hline
$N_p$ & 8 & 16 & 32 & 64 & 128 \\
\hline
$N_{timestep}$ & 249   & 248   & 239   & 246   & 243   \\
$N_{newton}$   & 2182  & 2178  & 2115  & 2151  & 2135  \\
$N_{gmres}$   & 64340 & 64567 & 64649 & 64039 & 63277 \\
$N_{newton}/N_{timestep}$ & 8.8 & 8.8 & 8.8 & 8.7 & 8.8 \\
$N_{gmres}/N_{newton}$ & 29.5 & 29.6 & 30.6 & 29.8 & 29.6 \\
\hline
\end{tabular}
\end{table}

\begin{figure}[!htbp]
  \centering
  \includegraphics[width=0.5\textwidth]{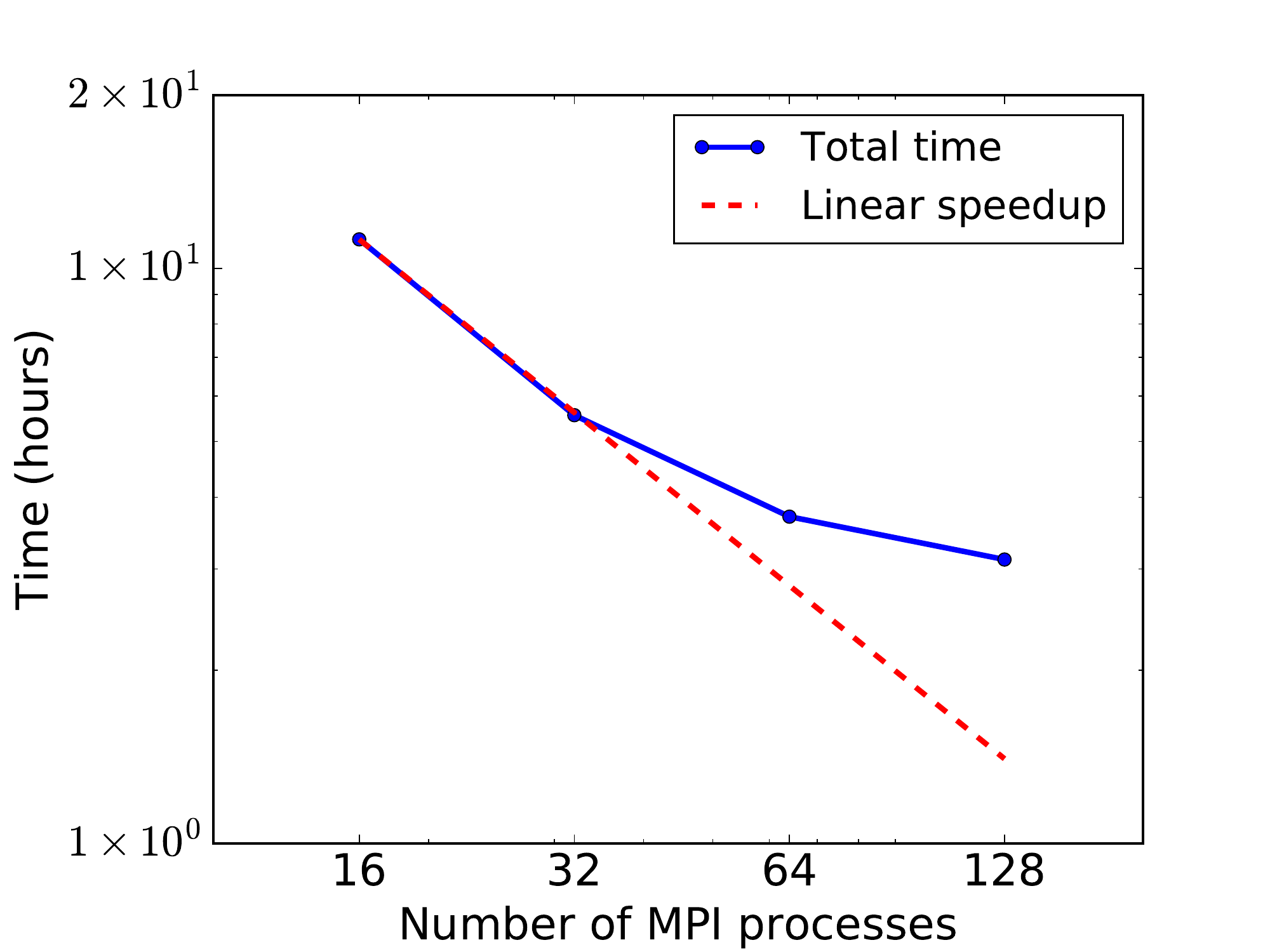}
  \caption{Total computation time vs. number of MPI processes for the black oil model test case with tetrahedral mesh.}
  \label{timetet}
\end{figure}

Figure~\ref{timetet} exhibits the total simulation times as a function of the number of MPI processes. 
The results are similar than in the previous test case. The scalability is very good up to $32$ MPI processes 
and degrades for $N_p=64$ and $128$ as expected for a number of unknowns in the pressure block roughly equal to $10^6$. {Table \ref{partitioning_tet} shows the maximum and mean number of own d.o.f. by process for $N_p=64,128$ with a larger disbalance for own nodes + fracture faces than in the previous test case but still quite smaller than the loss of parallel efficiency exhibited in Figure \ref{timetet} which is mainly due to the communication overhead.} 

\subsubsection{Water injection}

{
We modify the previous test case using the new fracture width $d_f= 0.001$ meter and injecting pure water instead of oil at the bottom boundary with a bottom pressure of $3$ MPa. 
The relative permeabilities are modified using a residual water saturation $S_{r}^w=0.2$ and
the initial water saturation is fixed to $S^w = S_{r}^w$.
}
{
The time stepping parameters are fixed to 
$t_f=5000$ days, 
$\Delta t^{(0)} = 0.001 \text{ days}$, $\Delta t^{(1)} = 30 \text{ days}$, 
$\Delta t^{(2)} = 100 \text{ days}$, $\Delta t^{(3)} = 30 \text{ days}$, 
$t^{(1)} = 600 \text{ days}$, $t^{(2)} = 2000 \text{ days}$, $k_f = 3$.
}

\begin{figure}[!htbp]
  \centering
  \includegraphics[width=0.4\textwidth]{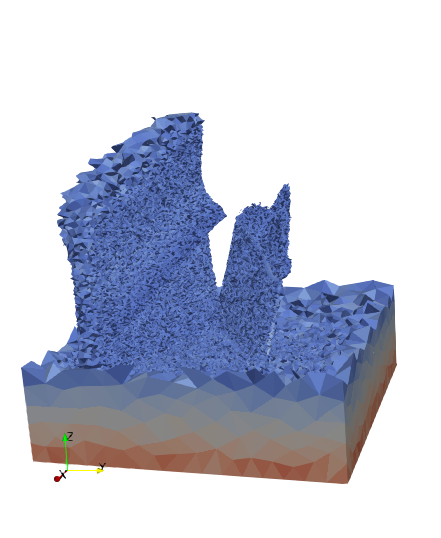}
  \hspace{0.1\textwidth}
  \includegraphics[width=0.4\textwidth]{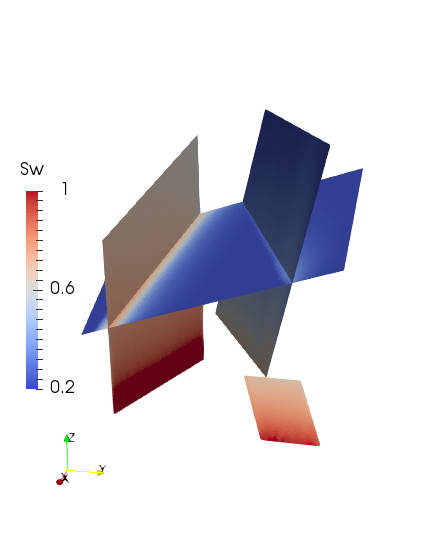}
  \caption{Water saturation at time $t_f$ in the matrix (left figure with a threshold of $S^w = 0.3$) and in the fractures (right) for the Black oil test case with water injection. }
  \label{BOWaterInjection}
\end{figure}

%

\begin{table}[!htbp]
\centering
\caption{Number of time steps ($N_{timestep}$), total number of Newton iterations ($N_{newton}$) and total number of linear solver iterations ($N_{gmres}$) vs. number of MPI processes / mesh size for the water injection black oil model test case.}
\label{Nitertet_weakscaling}
\begin{tabular}{c|c|c|c|c|c}
\hline
$N_p$ / nb of cells & 16/$1.17 \times 10^6$ & 32/$2.03 \times 10^6$ & 64/$4.11 \times 10^6$  \\
\hline
$N_{timestep}$            & 187   & 187  &  190    \\
$N_{newton}$              & 633   & 647  & 688  \\
$N_{gmres}$               &  8841  & 10168  & 12889  \\
CPU time  (s)            &  1213  &  1424 & 1940 \\
CPU time / $N_{gmres}$    &  0.137  & 0.140 & 0.151 \\
\hline
\end{tabular}
\end{table}

{
  In order to investigate the weak scalability of the code, Table \ref{Nitertet_weakscaling} exhibits the numerical behavior of the simulation obtained
  for this test case using tetrahedral meshes with $1.17 \times 10^6$, $2.03 \times 10^6$, $4.11 \times 10^6$ cells
  on respectively $N_p=16,32,64$ processes. 
  The number of Newton iterations as well as the total number of GMRES iterations increase only moderately with the mesh size.
  The CPU time per GMRES iteration exhibits a good weak
  scalability for $N_p = 16,32,64$ processes. 
}

\subsection{Non-isothermal liquid-gas simulation with a large discrete fault network}

We consider in this subsection a single $H_2O$ component liquid-gas non-isothermal model with $\mathcal{P}=\{\mathrm{water, gas}\}$ and $\mathcal{C}=\{ H_2O\}$. 
The thermodynamical laws providing the phase molar densities, viscosities, internal energies, 
and enthalpies as well as the saturation vapor pressure are obtained from \cite{refthermo}. The thermal conductivity is fixed to $\lambda = 2$ W m$^{-1}$ K$^{-1}$ 
and the rock volumetric internal energy is defined by 
$E_r(T) = c_p^r T$ with $c_p^r = 16. 10^{5}$ J m$^{-3}$ K$^{-1}$.  The gravity is not considered in this test case 
which means that the solution is 2 dimensional. 

The fault network is provided by M. Karimi-Fard and A. Lap\`ene from Stanford University and TOTAL as well as  
the prismatic mesh of the domain $\Omega=(0,5888)\times (0,3157) \times (0,200)$ (meters) which 
contains about $1.3\times 10^6$ prismatic cells, $3.4\times 10^6$ nodes and $7.1\times 10^5$ fault faces. 
The 3D mesh is defined by the tensor product of a triangular 2D mesh with a uniform vertical 1D mesh with $10$ intervals.  
The fault network contains $581$ connected components. The fault width is set to 
$d_f = 1$ m and the permeabilities are isotropic and fixed to $\Lambda_m=10^{-15} $ m$^2$ in the matrix domain and to 
$\Lambda_f=10^{-12} $ m$^2$ in the fault network. The porosities in the matrix domain and in the faults are $\phi_m=0.1$ and $\phi_f=0.1$ respectively.

\begin{figure}[!htbp]
  \centering
 \includegraphics[width=0.42\textwidth]{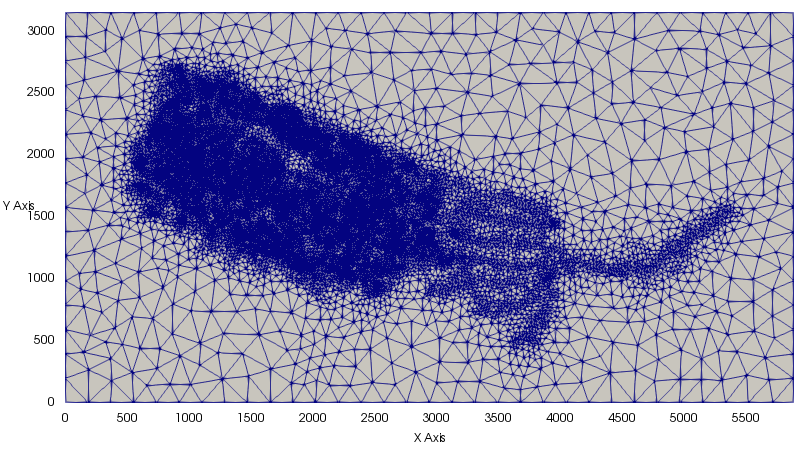}
  \hspace{0.1\textwidth}
  \includegraphics[width=0.42\textwidth]{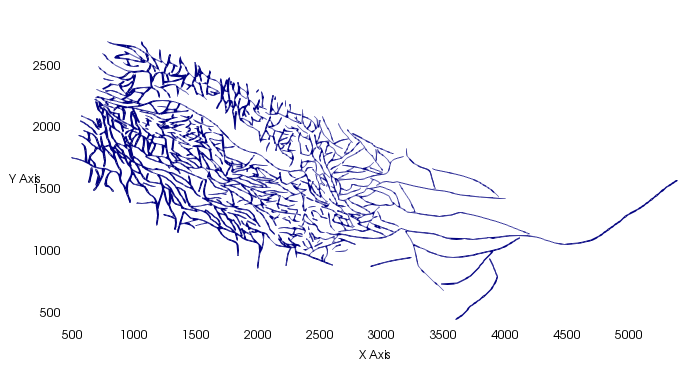}
  \caption{Horizontal view of the prismatic mesh of the matrix domain coarsened by a factor $2$ in the $x,y$ directions  (left) 
    conforming to the fault network (right).}
  \label{meshggf}
\end{figure}

\begin{table}[!htbp]
\centering
\caption{{Maximum/mean number of own cells, own cells+nodes+fracture faces and own nodes+fracture faces by process for the prismatic mesh with $1.3\times 10^6$ cells and $N_p=64,128$.} }
\label{partitioning_prism}
\begin{tabular}{c|c|c}
\hline
$N_p$ & 64 & 128  \\
\hline
own cells & 20214/20212 & 10107/10106    \\
own cells+nodes+fracture faces & 33767/33016 & 17010/16508   \\
own nodes+fracture faces & 13554/12803 & 6903/6401    \\
\hline
\end{tabular}
\end{table}

Let us set $\Gamma_{output}= \{(x,y,z) \in \Omega \ | \ x=0\}$ and $\Gamma_{input} = \{(x,y,z) \in \Omega \ | \ x=5888\}$.
The simulation domain is initially in liquid phase with $P=1$ MPa and $T=450$ K. Dirichlet boundary conditions are imposed 
at $\Gamma_{output}$ with $P=1$ MPa and $T=450$ K  (liquid phase) and at $\Gamma_{input}$ with $P=2$ MPa and $T=550$ K (gas phase). 
The remaining boundaries are considered impervious to mass and energy.

The linear and nonlinear solver parameters are fixed to $N_{newton}^{max}=50$, 
$N_{gmres}^{max}=300$, $\epsilon_{gmres}= 10^{-4}$, $\epsilon_{newton}= 10^{-6}$, 
and the time stepping parameters are fixed to $t_f=280000$ days, 
$\Delta t^{(0)} = 1000 \text{ days}$, $\Delta t^{(1)} = 10000 \text{ days}$, 
$k_f = 1$.

\begin{figure}[!htbp]
  \centering
  \addtocounter{subfigure}{2} 
  \begin{subfigure}{0.4\textwidth}
    \centering
    \includegraphics[width=\textwidth]{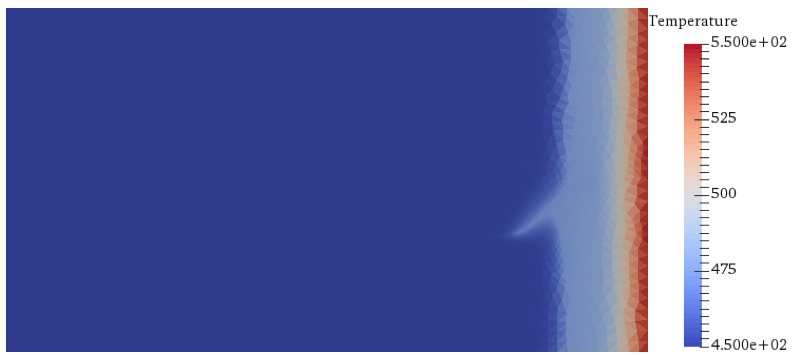}
    \caption{$t=70000$}
  \end{subfigure}
  \hspace{0.05\textwidth}
  \begin{subfigure}{0.4\textwidth}
    \centering
    \includegraphics[width=\textwidth]{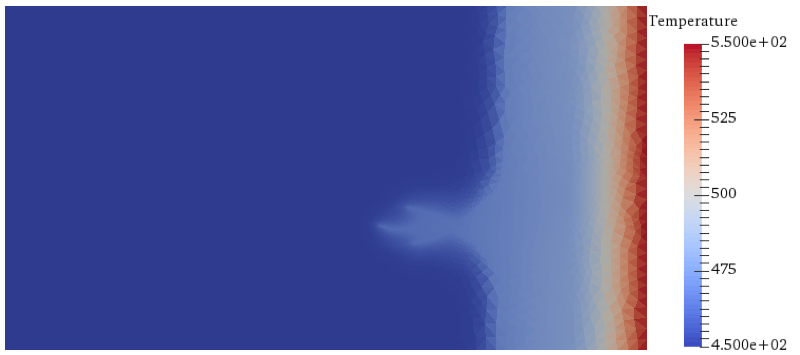}
    \caption{$t=140000$}
  \end{subfigure}
  \begin{subfigure}{0.4\textwidth}
    \centering
    \includegraphics[width=\textwidth]{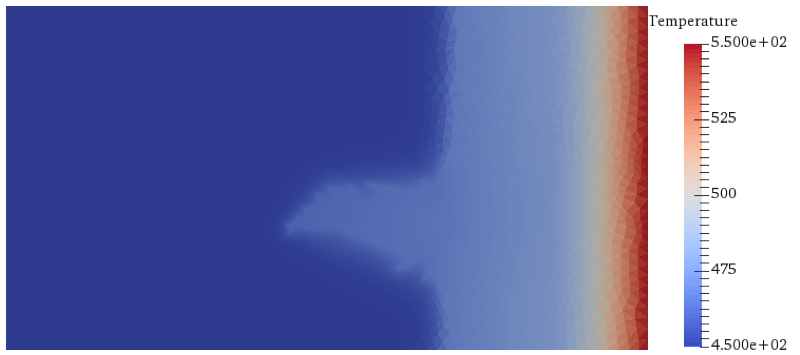}
    \caption{$t=210000$}
  \end{subfigure}
  \hspace{0.05\textwidth}
  \begin{subfigure}{0.4\textwidth}
    \centering
    \includegraphics[width=\textwidth]{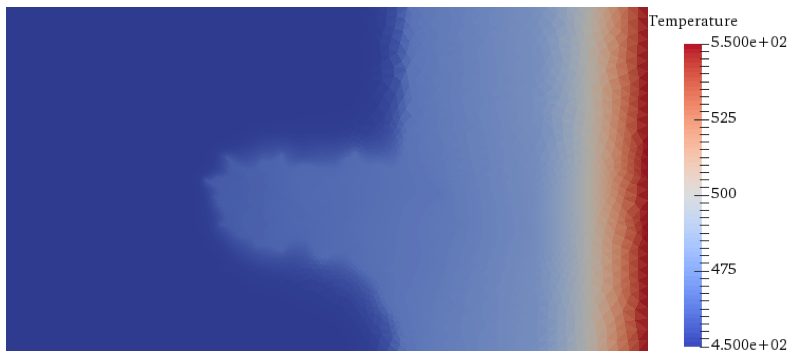}
    \caption{$t=280000$}
  \end{subfigure}
  \caption{Temperature in the matrix domain for the non-isothermal test case on the prismatic mesh.}
  \label{simuggfTemp}
\end{figure}

\begin{figure}[!htbp]
  \centering
  \addtocounter{subfigure}{2} 
  \begin{subfigure}{0.4\textwidth}
    \centering
    \includegraphics[width=\textwidth]{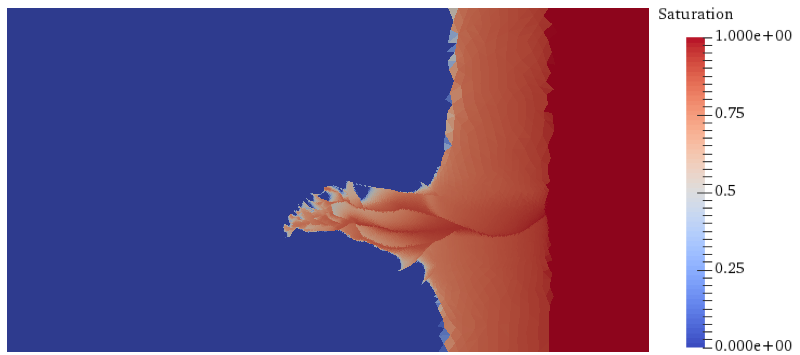}
    \caption{$t=70000$}
  \end{subfigure}
  \hspace{0.05\textwidth}
  \begin{subfigure}{0.4\textwidth}
    \centering
    \includegraphics[width=\textwidth]{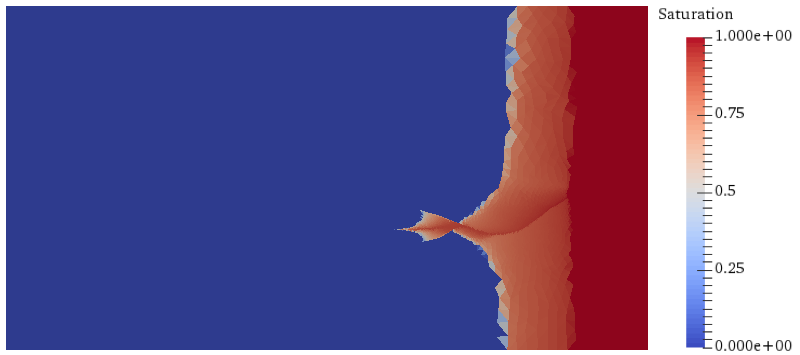}
    \caption{$t=140000$}
  \end{subfigure}
  \begin{subfigure}{0.4\textwidth}
    \centering
    \includegraphics[width=\textwidth]{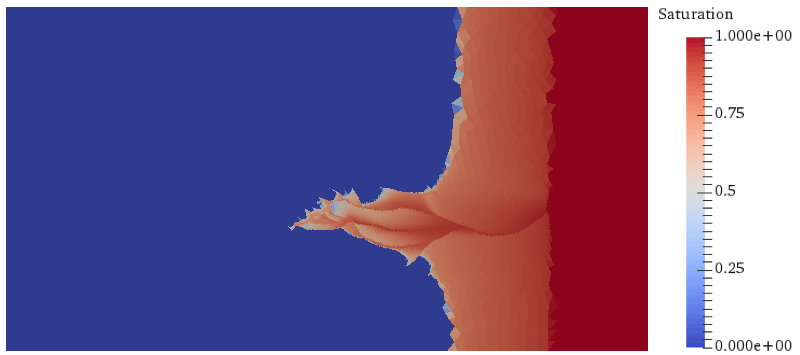}
    \caption{$t=210000$}
  \end{subfigure}
  \hspace{0.05\textwidth}
  \begin{subfigure}{0.4\textwidth}
    \centering
    \includegraphics[width=\textwidth]{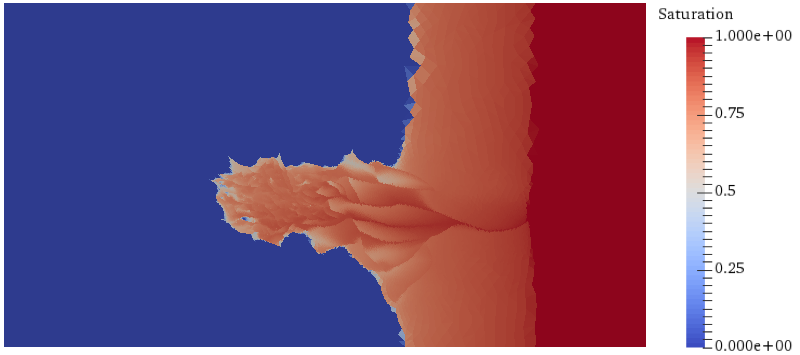}
    \caption{$t=280000$}
  \end{subfigure}
  \caption{Gas saturation in the matrix domain for the non-isothermal test case on the prismatic mesh.}
  \label{simuggfSat}
\end{figure}

\begin{figure}[!htbp]
  \centering
  \includegraphics[width=0.5\textwidth]{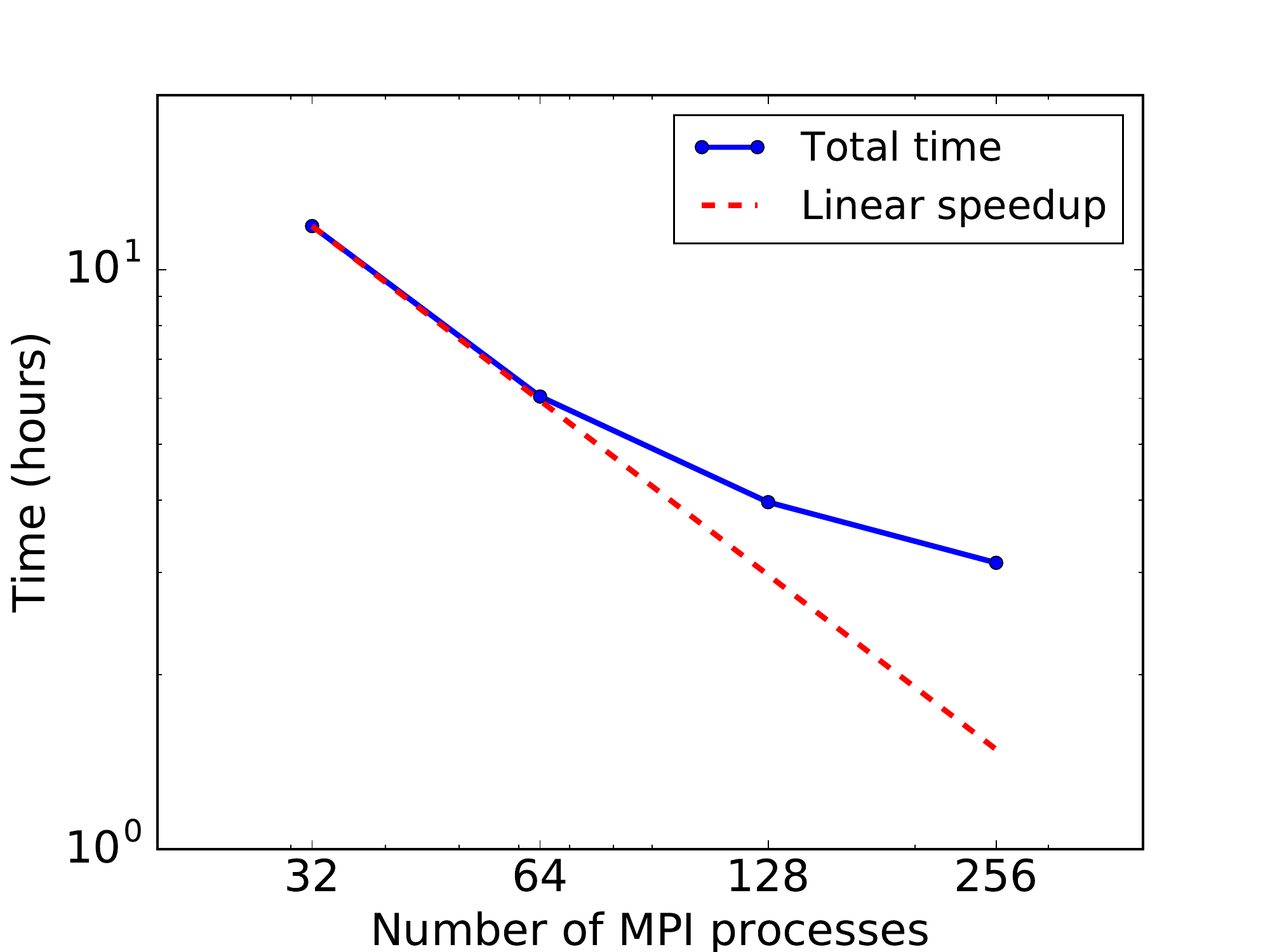}
  \caption{Total computation time vs. number of MPI processes for the non-isothermal test case on the prismatic mesh.}
  \label{timeggfmesh}
\end{figure}

Figures \ref{simuggfTemp} and \ref{simuggfSat} exhibit the temperature and the gas saturation at different times. 
Table \ref{Niterprism} shows the total number of Newton iterations 
and the total number of linear solver iterations which are, as for the previous test cases, almost 
independent on the number of MPI processes $N_p=32,64,128,256$. 
The average number of Newton iterations per time step is $20$. This is a high value but typical for 
such non-isothermal flows combining high non linearities in the thermodynamical laws 
and highly contrasted matrix and fault properties and scales. 
On the other hand, the number of linear solver iterations, roughly $60$ per Newton step, remains very good. 
Simarly as in the previous test cases, 
the scalability of the total simulation time with respect to the number of MPI processes presented in Figure \ref{timeggfmesh} 
is very good from $32$ to $64$ processes and then degrades for $N_p=128$ and $256$ due to a too small number of unknowns 
in the pressure block per MPI process.
{Table \ref{partitioning_prism} shows the maximum and mean number of own d.o.f. by process for $N_p=64,128$ with similar conclusions as in the previous test cases.} 

\begin{table}[!htbp]
\centering
\caption{Number of time steps ($N_{timestep}$), total number of Newton iterations ($N_{newton}$) and total number of linear solver iterations ($N_{gmres}$) vs. number of MPI processes for the non-isothermal test case on the prismatic mesh.}
\label{Niterprism}
\begin{tabular}{c|c|c|c|c}
\hline
$N_p$ & 32 & 64 & 128 & 256 \\
\hline
$N_{timestep}$ & 300    & 318    & 303    & 289 \\
$N_{newton}$   & 5890   & 6012   & 5946   & 5885 \\
$N_{gmres}$   & 372671 & 370954 & 383523 & 391250 \\
$N_{newton}/N_{timestep}$ & 19.6 & 18.9 & 19.6 & 20.4 \\
$N_{gmres}/N_{newton}$ & 63.3 & 61.7 & 64.5 & 66.5 \\
\hline
\end{tabular}
\end{table}

%

\subsection{Thermal convection test case with Cartesian mesh}

This test case considers the same physical two-phase non-isothermal model 
as in the previous subsection but including gravity. The simulation domain is $\Omega=(0,3000)^3$ in meters. The mesh is a 
3D uniform Cartesian mesh which contains $240^3$ cells. The fault network is defined by 
\begin{equation*}
\Gamma=\{(x,y,z)\in \Omega \ | \ x=1500 \ \mathrm{or} \ y=1500 
\ \mathrm{and} \  500 \leqslant z \leqslant 2500\}, 
\end{equation*}
with fault width fixed to $d_f=1$ meter. The permeabilities are isotropic and set 
to $\Lambda_m=10^{-15} $ m$^2$ in the matrix domain and to $\Lambda_f=10^{-12} $ m$^2$ in the fault network. The porosities in the matrix domain and in the faults are $\phi_m=0.1$ and $\phi_f=0.5$ respectively.

The domain is initially in liquid phase with a fixed temperature $293$ K 
and an hydrostatic pressure defined by its value $P=1$ bar at the top boundary. 
The temperature is fixed to $623$ K (liquid phase) at the bottom boundary which is impervious to mass. 
At the top boundary, the pressure is set to $1$ bar and the temperature to $293$ K (liquid phase).
A zero flux for both mass and temperature is imposed at the lateral boundaries of the domain. 

The linear and nonlinear solver parameters are fixed to $N_{newton}^{max}=25$, 
$N_{gmres}^{max}=300$, $\epsilon_{gmres}= 10^{-4}$, $\epsilon_{newton}= 10^{-5}$, 
and the time stepping parameters are fixed to $t_f=2\times 10^7$ days, 
$\Delta t^{(0)} = 5\times 10^5 \text{ days}$, $\Delta t^{(1)} = 5 \times 10^5  \text{ days}$, 
 $\Delta t^{(2)} = 10^5  \text{ days}$,  $\Delta t^{(3)} = 5\times 10^3  \text{ days}$, 
$t^{(1)} = 2.5 \times 10^6 \text{ days}$, $t^{(2)} = 1.7 \times 10^7 \text{ days}$,  $k_f = 3$. 

\begin{figure}[!htbp]
  \centering
  \addtocounter{subfigure}{2} 
  \begin{subfigure}{0.75\textwidth}
    \centering
    \includegraphics[width=\textwidth]{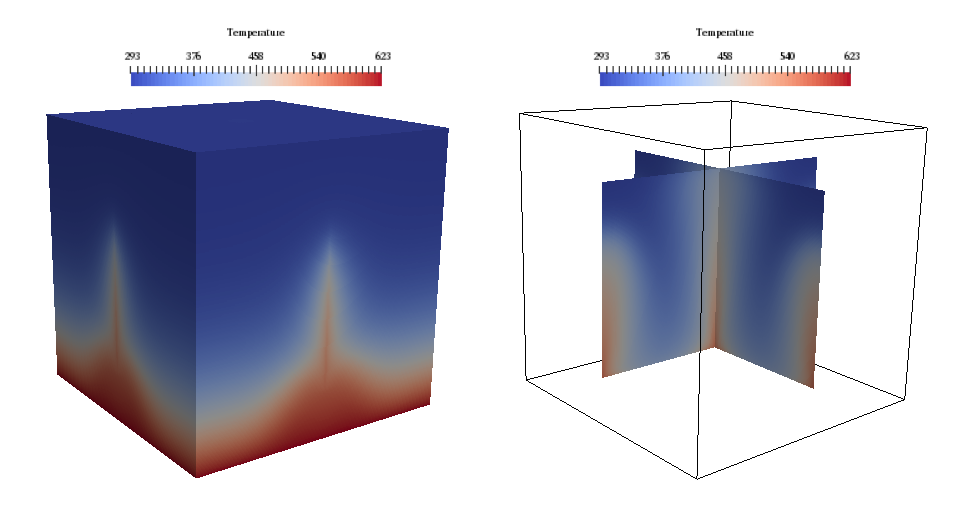}
    \caption{$t=1\times 10^7$ days}
  \end{subfigure}
  \hspace{0.0\textwidth}

  \begin{subfigure}{0.75\textwidth}
    \centering
    \includegraphics[width=\textwidth]{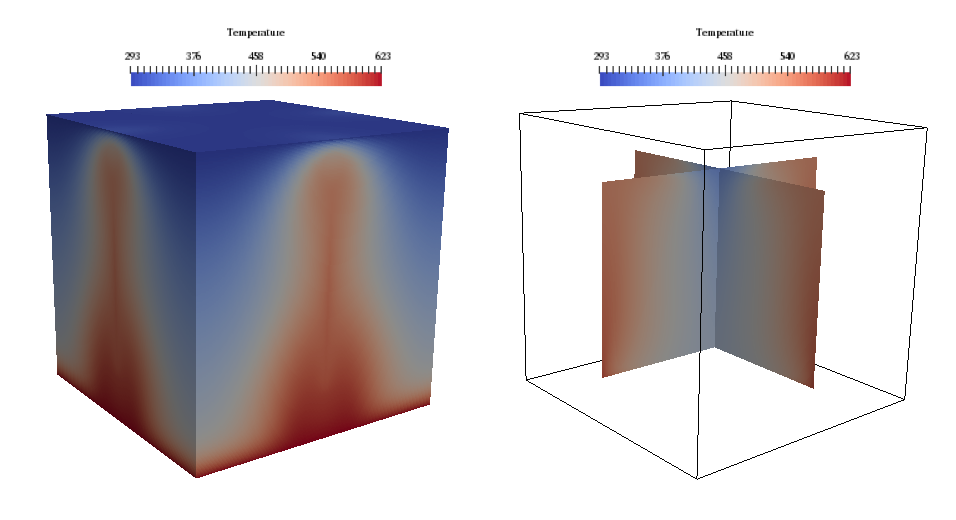}
    \caption{$t=2\times 10^7$ days}
  \end{subfigure}

  \caption{Temperature in the faults and in the matrix domain 
    at different times (days).}
  \label{simucar_Temp}
\end{figure}

\begin{figure}
  \centering
  \includegraphics[width=0.4\textwidth]{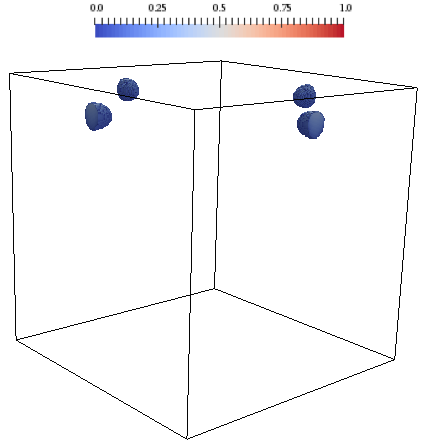}
  \caption{Saturation of gas in the matrix domain at the final time $t_f=2\times 10^7$ days.}
  \label{simucar_Sat}
\end{figure}

Figure \ref{simucar_Temp} shows the temperature in the faults and in the matrix domain at times $t=1\times 10^7$ {days} and $t=t_f$. In addition, we present in Figure \ref{simucar_Sat} the gas saturation at final time.

In this test case, the thermal convection leads to convective instabilities which are triggered by the numerical round-off errors. 
Hence it is not appropriate to make scalability tests since the solution will depend on 
the number of MPI processes. Therefore, we only exhibit in Table \ref{Niterunstable} the results obtained 
for $N_p=256$. This simulation lasts $40.57$ hours. The convective instabilities and strong nonlinearities require a small time step 
in order to obtain Newton's convergence, especially at the end of the simulation when the gas phase appears. 
\begin{table}[!htbp]
\newcolumntype{C}[1]{>{\centering\let\newline\\\arraybackslash\hspace{0pt}}m{#1}}
\centering
\caption{Number of time steps ($N_{timestep}$), total number of Newton iterations ($N_{newton}$) and total number of linear solver iterations ($N_{gmres}$) for the thermal convection test case on the Cartesian mesh where $N_p=256$.}
\label{Niterunstable}
\begin{tabular}{c|C{2cm}}
\hline
$N_{timestep}$ &  3117  \\
$N_{newton}$   & 6712   \\
$N_{gmres}$   &  238600 \\
$N_{newton}/N_{timestep}$ & 2.2 \\
$N_{gmres}/N_{newton}$ & 35.5 \\
\hline
\end{tabular}
\end{table}

\subsection{Thermal convection test case with tetrahedral mesh}

This last test case considers the same physical two-phase non-isothermal model as in the previous subsection,
but with the tetrahedral mesh shown in Figure \ref{meshtet} and rescaled to a larger domain $\Omega=(0,3000)^3$ in meters. 
The fault width is fixed to $d_f=1$ meter. The permeabilities are isotropic and set 
to $\Lambda_m=10^{-14} $ m$^2$ in the matrix domain and to $\Lambda_f=10^{-12} $ m$^2$ in the fault network.
The porosity in the matrix domain is set to $\phi_m=0.25$, and to $\phi_f=0.35$ in the fault network.

At the intersection $\{z=0\} \cap \Gamma$ of the bottom boundary with the fault network,
the temperature is fixed to $623$ K
and a mass flow rate of $100$ kg/s is uniformly prescribed. 
At the matrix bottom boundary, the temperature is set as $473$ K and the mass flow rate is set to zero.  
At the top boundary, the pressure is set to $10^5$ Pa and the temperature to $293$ K (liquid phase).
A zero flux for both mass and temperature is imposed at the lateral boundaries of the domain. 
The simulation domain is initially in liquid phase with an hydrostatic pressure defined by the pressure
boundary condition at the top boundary and with 
a linear temperature between $293$K at the top boundary and $473$K at the bottom boundary.

We set the linear and nonlinear solver parameters to $N_{newton}^{max}=30$, 
$N_{gmres}^{max}=150$, $\epsilon_{gmres}= 10^{-5}$, $\epsilon_{newton}= 10^{-6}$, 
and the time stepping parameters are fixed to $t_f=2\times 10^5$ days, 
$\Delta t^{(0)} = 50 \text{ days}$, $\Delta t^{(1)} = 1000  \text{ days}$, 
 $\Delta t^{(2)} = 50  \text{ days}$, 
$t^{(1)} = 5.8 \times 10^4 \text{ days}$,  $k_f = 2$.

\begin{figure}[!htbp]
  \centering
  \addtocounter{subfigure}{2} 
  \begin{subfigure}{0.75\textwidth}
    \centering
    \includegraphics[width=\textwidth]{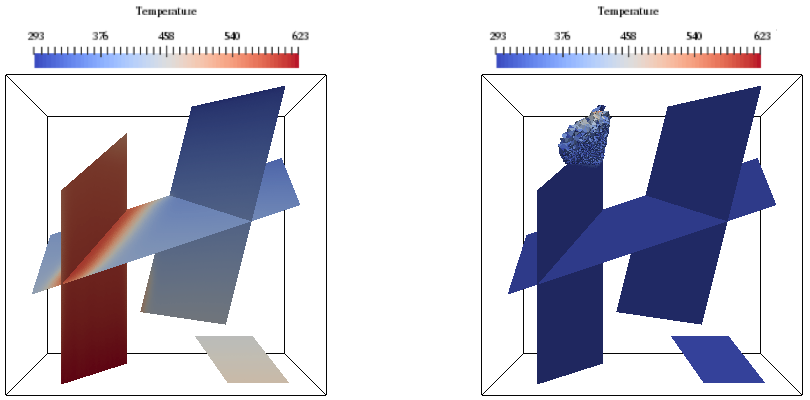}
    \caption{$t=7\times 10^4$ days}
  \end{subfigure}
  \hspace{0.0\textwidth}

  \begin{subfigure}{0.75\textwidth}
    \centering
    \includegraphics[width=\textwidth]{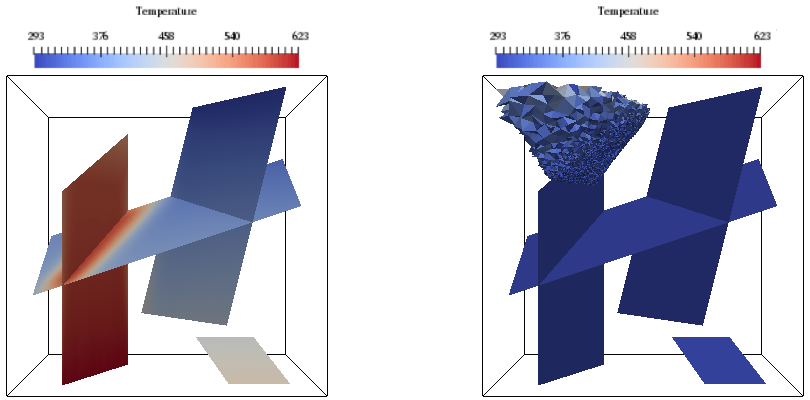}
    \caption{$t=2\times 10^5$ days}
  \end{subfigure}

  \caption{Left: temperature in the faults at different times (days). Right: saturation of gas in the faults and in the matrix domain at different times (days) where a threshold of 0.0001 is used on the matrix domain.}
  \label{simutet_Temp}
\end{figure}

Figure \ref{simutet_Temp} exhibits the temperature in the faults (left) and the gas saturation 
in the faults and in the matrix domain (right) at times $t=7\times 10^4$ days and $t=t_f$. 

\begin{figure}[!htbp]
  \centering
  \includegraphics[width=0.5\textwidth]{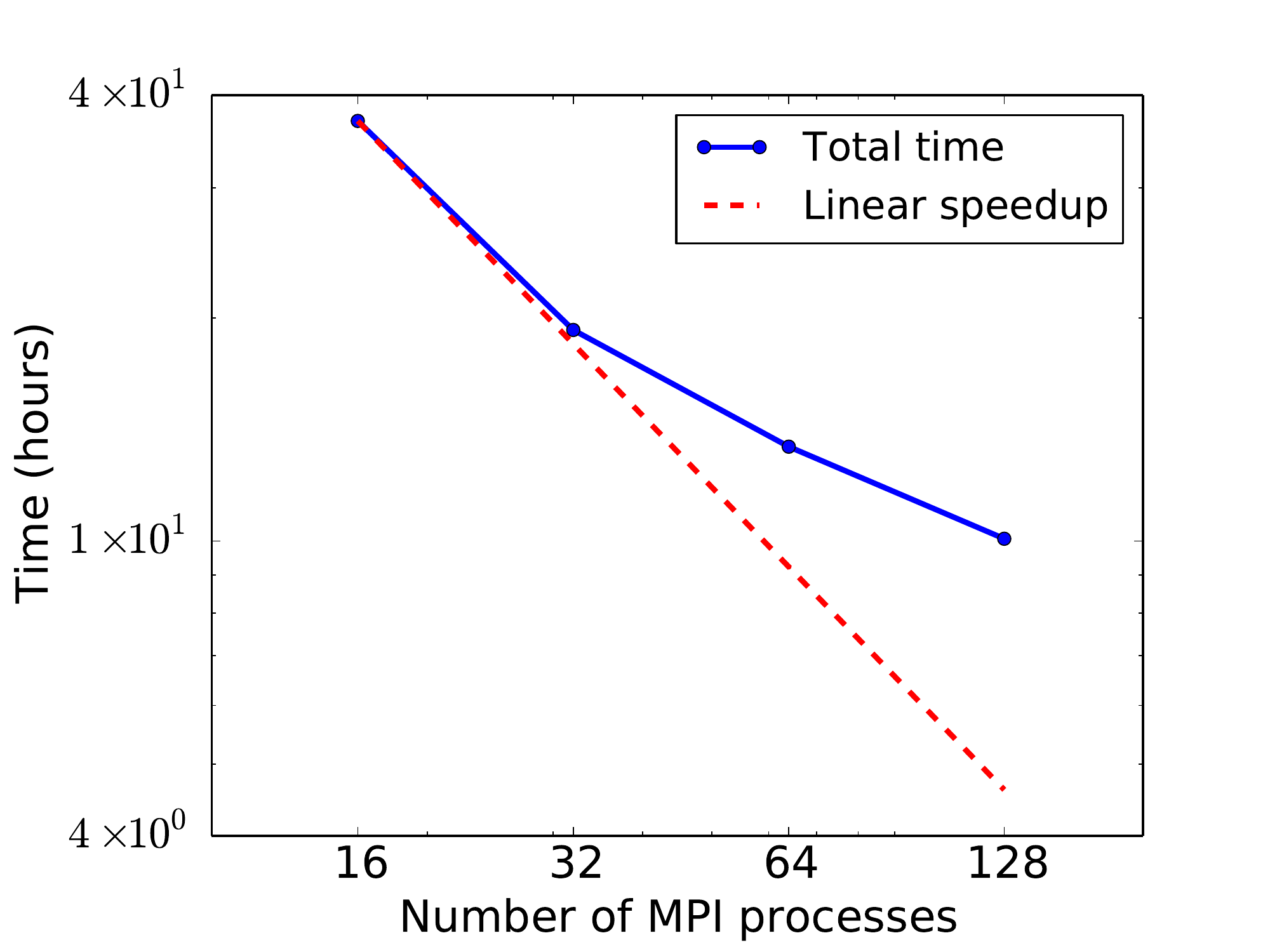}
  \caption{Total computation time vs. number of MPI processes for the thermal convection test with tetrahedral mesh.}
  \label{timetemptet}
\end{figure}

Table \ref{Nitertettemp} shows the total number of Newton iterations 
and the total number of linear solver iterations. We present the total computation time in hours for different number of MPI processes $N_p=16,32,64,128$ in Figure \ref{timetemptet}.
The scalability is similar to the one obtained with
the black oil model test case using the same tetrahedral mesh as shown in Figure \ref{timetet}.
\begin{table}[!htbp]
\centering
\caption{Number of time steps ($N_{timestep}$), total number of Newton iterations ($N_{newton}$) and total number of linear solver iterations ($N_{gmres}$) vs. number of MPI processes for the thermal convection test case with tetrahedral mesh.}
\label{Nitertettemp}
\begin{tabular}{c|c|c|c|c}
\hline
$N_p$ & 16 & 32 & 64 & 128  \\
\hline
$N_{timestep}$ & 2939 & 2946 & 2950 & 2946 \\
$N_{newton}$   & 10043 & 10259 & 10352 & 10074 \\
$N_{gmres}$    & 151993 & 155925 & 158124 & 153836 \\
$N_{newton}/N_{timestep}$ & 3.4 & 3.5 & 3.5 & 3.4 \\
$N_{solver}/N_{newton}$ & 15.1 & 15.2 & 15.3 & 15.3 \\
\hline
\end{tabular}
\end{table}

\section{Conclusion}

In this paper, a discrete fracture model accounting for non-isothermal compositional multiphase  
Darcy flows was introduced. The geometry takes into account complex networks of intersecting, immersed or non immersed 
planar fractures. The physical model accounts for an arbitrary nonzero number of components in each phase allowing to 
model immiscible, partially miscible or fully miscible flows. The discretization is based on the VAG finite volume scheme adapted 
to unstructured polyhedral meshes and to anisotropic heterogeneous media. The time integration is fully implicit in order to avoid strong restrictions on the time step due to the high velocities and small volumes in the fractures. The discrete model is implemented in parallel based on 
the SPMD paradigm and using one layer of ghost cells in order to assemble the systems locally on each processor. The CPR-AMG preconditioner was investigated to deal with non-isothermal models. 

The numerical results exhibit the ability of our  discrete model 
to combine complex physics including non-isothermal flows, thermodynamical equilibrium and buoyancy forces  
with fracture networks including highly contrasted matrix fracture permeabilities. 
The parallel scalability requires, as expected for fully implicit discretizations when using AMG type preconditioners, 
that the number of degrees of freedom per processor is kept high enough.

\section*{Acknowledgments} 

This work was supported by a joint project between INRIA and BRGM Carnot institutes (ANR, INRIA, BRGM) and
partially supported by the CHARMS ANR project (ANR-16-CE06-0009). This work was also granted access to the HPC and visualization resources of ``Centre de Calcul Interactif" hosted by University Nice Sophia-Antipolis. 

\bibliographystyle{elsarticle-num}
\bibliography{Bib-revised}
\end{document}